\documentclass[11pt]{amsart}
\usepackage[T1]{fontenc}
\usepackage[utf8]{inputenc}
\usepackage{amssymb,amsmath}
\usepackage{esint} \usepackage{graphicx}
\usepackage{bm,bbm}
\usepackage{color}
\usepackage{tikz}
\usepackage{float}
\usepackage{framed}
\usepackage{makecell}  \usetikzlibrary{snakes}
\usepackage{pgfplots}
\pgfplotsset{compat=newest}
\usepackage{array,blkarray}
\usepackage{bigdelim}
\usepackage{listings} 
\usepackage{hhline}
\usepackage{arydshln}
\usepackage{multirow}
\usepackage{enumitem}

\usepackage[hidelinks]{hyperref}

\usepackage{mathtools}  \newtheorem{theorem}{Theorem}

\theoremstyle{definition}
\newtheorem{definition}[theorem]{Definition}

\theoremstyle{remark}
\newtheorem{remark}[theorem]{Remark}

\numberwithin{theorem}{section}
\numberwithin{equation}{section}
\numberwithin{table}{section}
\numberwithin{figure}{section}

\def\R{\mathbb{R}}

\def\qo{\ensuremath{v_0}}

\def\eps{\varepsilon}

\definecolor{color0}{rgb}{0.17004232121057961, 0.43679759647517286, 0.22372555555555548}
\definecolor{color1}{rgb}{0.63284224750184226, 0.4747981096220677, 0.29070209208025455}
\definecolor{color2}{rgb}{0.82995767878942039, 0.56320240352482709, 0.77627444444444449}

\definecolor{blau}{RGB}{0, 51, 255}
\definecolor{hellblau}{RGB}{153, 204, 255}
\definecolor{hellrot}{RGB}{255, 0, 0}
\definecolor{rot}{RGB}{156, 0, 0}          \definecolor{firebrick}{RGB}{176, 34, 34} 
\definecolor{deep_pink}{RGB}{255, 20, 147} 
\definecolor{sky_blue}{RGB}{74, 112, 139}
\definecolor{slate_blue}{RGB}{71, 60, 139}
\definecolor{chartreuse}{RGB}{118, 238, 0}
\definecolor{chartreuseL}{RGB}{228, 255, 150}
\definecolor{light_blue}{RGB}{178, 223, 238}
\definecolor{dodge_blue}{RGB}{17, 78, 138}
\definecolor{code_backg}{RGB}{238, 216, 174}

\definecolor{myBlue1}{RGB}{101,149,239}  \definecolor{myBlue2}{RGB}{113,104,238} \definecolor{myBlue3}{RGB}{30,144,255} \definecolor{myGreen1}{RGB}{154,204,50} \definecolor{myGreen2}{RGB}{69,169,0} \definecolor{myGreen3}{RGB}{154,205,50} \definecolor{myGreen4}{RGB}{105,139,34} \definecolor{myRed1}{RGB}{210,105,30} \definecolor{myRed2}{RGB}{165,42,42} \definecolor{myRed3}{RGB}{139,26,26} \definecolor{myLGray}{RGB}{225,225,225}

\DeclareMathOperator{\diag}{diag}
\DeclareMathOperator{\ddiv}{div}

\DeclareMathOperator{\rank}{rank}

\DeclareMathOperator{\CR}{CR}

\newcommand{\dx}{\ensuremath{\ \text{d}x }}

\newcommand{\calB}{\ensuremath{\mathcal{B}} }

\newcommand{\calE}{\ensuremath{\mathcal{E}} }

\newcommand{\calF}{\ensuremath{\mathcal{F}} }
\newcommand{\calG}{\ensuremath{\mathcal{G}} }
\newcommand{\cH}{\ensuremath{\mathcal{H}} }

\newcommand{\calK}{\ensuremath{\mathcal{K}} }

\newcommand{\calP}{\ensuremath{\mathcal{P}} }

\newcommand{\Q}{\ensuremath{\mathcal{Q}} }

\newcommand{\calT}{\ensuremath{\mathcal{T}} }
\newcommand{\V}{\ensuremath{\mathcal{V}} }

\newcommand{\calW}{\ensuremath{\mathcal{W}} }

  \newcommand{\inicond}{\ensuremath{a}}

\newcommand{\Gammai}{\Gamma_\text{i}}
\newcommand{\Gammao}{\Gamma_\text{o}}
\newcommand{\Gammaw}{\Gamma_\text{w}}

\def\ds{\,\text{d}s}

\def\dx{\,\text{d}x}

\newcommand{\ddt}{\ensuremath{\frac{\text{d}}{\text{d}t}} }
 \DeclareMathOperator{\dive}{div}
\def\trp{{\mathsf T}}
\def\grad{\nabla}

\providecommand{\inva}[1]{\text{~\textup{d}} #1}
\DeclarePairedDelimiter\norm{\|}{\|}

\providecommand{\pdt}[1]{\frac{\partial #1}{\partial t}}

\def\RE{\text{Re}}
\def\v{\mathfrak{v}}
\def\p{\mathfrak{p}}
\def\rhs{\text{rhs}}
\def\vS{v_\text{S}}
\def\pS{p_\text{S}}
\def\pNS{p_\text{NS}}

\def\EApair{\ensuremath{\bigl ( \mathcal E, \mathcal A \bigr )}}
\def\dDAE{\ensuremath{\Delta}AE}

\def\lintol{\texttt{tol}} 

\def\ver{e_{\tau;\lintol}^v}
\def\per{e_{\tau;\lintol}^p}
\def\resmom{r_{\tau;\lintol}^M}
\def\rescon{r_{\tau;\lintol}^C}
\def\rfrc{\texttt{ref}}
\DeclareMathOperator{\trapz}{trp}

\definecolor{color0}{rgb}{0.0964154061637, 0.177286718596, 0.283473631854}
\definecolor{color1}{rgb}{0.170042321211, 0.436797596475, 0.223725555556}
\definecolor{color2}{rgb}{0.632842247502, 0.474798109622, 0.29070209208}
\definecolor{color3}{rgb}{0.829957678789, 0.563202403525, 0.776274444444}
\definecolor{color4}{rgb}{0.763267446284, 0.850242277575, 0.951553976205}
 
\begin{document}
\title[Continuous and Discrete Navier-Stokes Equations]{Continuous, Semi-discrete, and Fully Discretized Navier-Stokes Equations}
\author[]{R.~Altmann$^{*}$ and J.~Heiland$^{\dagger}$}
\address{${}^{*}$ Institut f\"ur Mathematik MA4-5, Technische Universit\"at Berlin, Stra\ss e des 17.~Juni 136, 10623 Berlin, Germany}
\address{${}^{\dagger}$ Max Planck Institute for Dynamics of Complex Technical Systems, Sandtorstra{\ss}e 1, 39106 Magdeburg, Germany}
\email{raltmann@math.tu-berlin.de, heiland@mpi-magdeburg.mpg.de}
\date{\today}
\begin{abstract}
The Navier--Stokes equations are commonly used to model and to simulate flow phenomena. We introduce the basic equations and discuss the standard methods for the spatial and temporal discretization. We analyse the semi-discrete equations -- a semi-explicit nonlinear DAE -- in terms of the strangeness index and quantify the numerical difficulties in the fully discrete schemes, that are induced by the strangeness of the system. By analyzing the Kronecker index of the difference-algebraic equations, that represent commonly and successfully used time stepping schemes for the Navier--Stokes equations, we show that those time-integration schemes factually remove the strangeness. The theoretical considerations are backed and illustrated by numerical examples. 
\end{abstract}
\maketitle

\section{Introduction}
The \emph{Navier--Stokes equations} (NSE) are a system of nonlinear partial-differential equations that have been commonly used to model fluid flows for more than a century. The NSE are believed to describe all kinds of incompressible flows sufficiently well as long as the setup supports the hypothesis that the fluid is a continuum. Indeed, comparisons of numerical simulations with experiments show arbitrarily good agreement of the model with the observations over a long range from slowly moving flows in small geometries like a pipe up to highly turbulent flows over large spatial extensions like the flow around an airplane or even weather phenomena. Nevertheless, the mere equations and the extent of their applicability have not been fully deciphered by now and a substantial progress in this respect will be eligible for a Clay price\footnote{see \url{http://www.claymath.org/millennium-problems/navier–stokes-equation}}.

Under the assumption of continuity of the observed quantities, the NSE can be derived from fundamental laws of physics; see~\cite{LanL87} and~\cite{ChoM93}. These considerations are well backed for a macroscopic viewpoint, from which a fluid like water appears as a continuum. On a microscopic level, where discrete molecular structures define the states, the NSE cannot capture the physics right, as it is well-known, e.g., for capillary flows. 

On the molecular level, fluids are better described by the \emph{Boltzmann equations}, which model molecular interactions. This fact seems undisputed the more that the NSE can also be interpreted and derived through a limiting process of the Boltzmann equations in the sense of averaging the microscopic quantities for a macroscopic description \cite{Sai09}.

As a mathematical object the NSE have ever been subject to fundamental investigations and led to its own research field and its own subject definitions in the \emph{MSC} classification scheme\footnote{see \url{http://www.ams.org/mathscinet/msc/msc.html?t=35Q30&btn=Current} and related }. The research on the NSE has focussed on the analysis of the equations and their numerical approximation. Early results on the existence of solutions are due to Leray \cite{Ler33} (weak solutions) and Fujita\&Kato \cite{FujK64} (smooth solutions). The first textbooks on the functional and on the numerical analysis were written by Ladyzhenskaya \cite{Lad69} and Temam \cite{Tem77}, respectively. 

For the numerical analysis of the spatial discretization, one may distinguish two lines of development. The mathematical line focusses on Galerkin methods in the realm of variational formulations whereas the engineering orientated line has been advancing finite volume methods (FVM) as they appear well-suited for simulations. On the side of Galerkin methods, and in particular finite element methods (FEM), there have been many efforts in designing stable elements like the famous \emph{Taylor--Hood} elements~\cite{TayH73} as well as for general convergence results; see the textbooks~\cite{GirR86,Lad69,Pir89} for the numerical analysis and \cite{Tur99} for an application oriented overview. On the side of FVM that are the method of choice in most general purpose flow solvers, there have been general developments in view of discretizing conservation laws \cite{LeV92} and particular progress in view of stable approximation of fluid flow \cite{FerP02}. 

A numerical analysis of approximations to the time-dependent NSE with FEM semi-discretizations has been carried out by Heywood\&Rannacher~\cite{HeyR90}; see also the textbook~\cite{GreS00}, that covers implementation issues. For time marching schemes for FVM formulations we refer to~\cite{FerP02}. Strategies for the iterative solution of the arising linear systems can be found in~\cite{ElmSW05} (FEM) or~\cite{FerP02} (FVM).

In this work we revisit the NSE from a differential-algebraic equations (DAE) perspective. This includes the modelling where the incompressibility is treated as an algebraic constraint in an abstract space and the spatial semi-discretization, that has to be handled with care to respect the incompressibility constraint and to lead to a well-posed classical DAE. It also includes the temporal discretization, in which the DAE properties of the NSE become evident and (hopefully not) problematic. Such a pure DAE perspective has been taken on by Weickert \cite{Wei97} who analysed finite difference approximations and certain time-stepping schemes for the NSE, by Emmrich\&Mehrmann \cite{EmmM13} who provided conditions and solution representations for linearized NSE in abstract spaces in line with linear time-invariant DAEs in finite-dimensional state-spaces. In \cite{Hei14} the nonlinear NSE and its Galerkin approximations have been analysed in view of consistency of reformulations of semi-discrete DAE approximations with the infinite-dimensional model. 

The paper is organized as follows. 
In Section~\ref{sect_pde} we derive the NSE from first principles and formulate the weak form as an {\em operator DAE}. 
The direct connection to DAEs is then made in Section~\ref{sect_dae}, in which we report on several spatial discretization schemes and that commonly used FEM schemes lead to systems of {\em strangeness index} one. 
In Section~\ref{sect_solution} we analyse the time approximation schemes in terms of the index of the resulting \emph{difference-algebraic equations} (\dDAE). We show that the straight-forward temporal discretization leads to a scheme of higher index than that of well-established time-stepping schemes for incompressible flows. We further confirm, that schemes with a \dDAE~of lower index can also be obtained from a standard time-discretization applied to a reformulation of the DAE with lower index.
In the numerical examples in Section~\ref{sect_numerics}, we confirm the superiority of the lower index \dDAE~approximations of the NSE over the straight-forward time-discretization. As a benchmark for time-integration schemes for the considered class of nonlinear semi-explicit DAEs of strangeness index 1, we provide reference trajectories and the system coefficients and nonlinear inhomogeneities for direct realization in \emph{Python}. We conclude this paper with summarizing remarks in Section~\ref{sect_conclusion}.
 \section{Continuous model}\label{sect_pde} 
\subsection{Derivation of the Navier--Stokes equations}\label{sect_NSE_derivation}
The Navier--Stokes equations (NSE) provide a model of a flow as a continuum. The basic assumption for their derivation and, thus, the validity of the model is that the flow under consideration forms a continuous entity of flow particles in a spatial domain $\Omega\subset \mathbb R^{3} $ and a time interval $\mathcal I$ such that the functions 
\begin{equation}
  \v \colon \mathcal I \times \Omega \to \mathbb R^{3}, \quad \p \colon \mathcal I \times \Omega \to \mathbb R, \quad\text{and}\quad \rho \colon \mathcal I \times \Omega \to \mathbb R ,
\end{equation}
describing the velocity, the pressure, and the density as measured at the position $x\in \Omega$ at time $t\in\mathcal I$ are continuous functions. Here, continuous means that we can apply differential calculus in order to derive basic partial differential equations. Later, when we derive the weak formulation of the NSE, the needed continuity will be specified further.

As mentioned in \cite[p.~2]{ChoM93}, these assumptions lead to a model that is believed provide accurate descriptions of common macroscopic flow phenomena. In, e.g., setups of small geometric scales like in capillary flows or under vacuum-like conditions \cite{KarBN05}, the discrete microscopic molecular structure of the fluid that constitutes the flow has to be taken into account.

The basic assumption of continuity of the matter allows for the consideration of a possibly infinitesimal small control volume $W\subset \Omega$, an open bounded domain in $\mathbb R^{3}$ that contains a given agglomerate of fluid particles in the considered flow. Continuity also implies that a fixed $W$ is deformed and convected by the flow but always consists of the same fluid particles.

Under this continuity assumption, one can call on the \emph{Reynolds Transport Theorem}, that relates the temporal change of an integral quantity over $W$ to the convection velocity $\v$. 
\begin{theorem}[see~\cite{Rey03}, Eq.(16) and \cite{ChoM93}, p.10]
  Let $W\colon I \to \mathbb R^{3}$ describe a smoothly moving control volume and let $f
  \colon t \times W(t) \mapsto \mathbb R^{}$ be a sufficiently smooth function, then
  \begin{equation}\label{eqn_RTT}
	\ddt \int_W f \inva V = \int_W \frac{\partial f}{\partial t} + \dive (f \v) \inva V.
  \end{equation}
\end{theorem}
\subsubsection*{Incompressibility and mass conservation}
The \emph{Reynolds Transport Theorem} can be used to show that a flow is \emph{incompressible}, which means that the volume of any agglomerate $W$ is constant over time, if and only if the velocity field $v$ is divergence free. In fact, if \eqref{eqn_RTT} applies, then one has that
\begin{equation}\label{eqn_flowdivfree}
  \ddt \int_W\inva V = 0 \quad\text{if, and only if,}\quad \int_W \dive \v \inva V=0.
\end{equation}
With a well-defined density function $\rho$, the mass of a control volume is defined as the integral over the (mass) density $\rho$ and, with the assumption that in the flow there are neither mass sinks nor mass sources, an application of \eqref{eqn_RTT} gives
\begin{equation}\label{eqn_contieqint}
  \ddt \int_{W} \rho \inva V =0 \quad\text{if, and only if,}\quad \int_W \frac{\partial \rho}{\partial t} + \dive (\rho \v) \inva V = 0.
\end{equation}
\begin{remark}
The conservation of mass and the incompressibility of a flow are closely related but only equivalent in the case that the density $\rho$ is constant in space and time. Flow models that assume incompressibility and varying density functions $\rho$ are applied, e.g., in oceanography~\cite{Tar06}.
\end{remark}
\begin{remark}
	If the volume of a fluid parcel changes over time, the flow is called \emph{compressible}. In this case, the mass density $\rho$ is modelled as an unknown function and related to the pressure $p$ through constitutive or so-called \emph{state equations}, like the \emph{ideal gas low}; see \cite{FeiKP16}.
\end{remark}
\subsubsection*{Balance of momentum}
Under the continuum assumption, the momentum of a fluid agglomerate $W$ can be expressed as the integral of the mass density times velocity and the temporal change equated with volume and surface forces on $W$:
\begin{equation}\label{eqn_mombalance}
  \ddt \int_{W} \rho \v \inva V = \int_W \rho g\inva V + \int_{\partial W} \sigma n \inva S, 
\end{equation}
where $g\colon \mathcal I \times \Omega \to \mathbb R^{3}$ is the density of a body force, and where $\sigma \colon \mathcal I \times \Omega \to \mathbb R^{3,3}$ is a tensor such that $\sigma n$, where $n$ is the normal field on the boundary $\partial W$ of $W$, represents the density of the forces acting on the surface. Note that the expression in \eqref{eqn_mombalance} is vector valued and that integration and differentiation is performed componentwise.

Applying the \emph{Divergence Theorem} componentwise, one can write the term with the surface forces as a volume integral
\begin{equation}\label{eqn_divthmstresstens}
  \int_{\partial W} \sigma n \inva S= \int_{W} \dive \sigma \inva V, 
\end{equation}
with the divergence of a tensor defined accordingly.

So far, all assumptions have based on first principles. For the mathematical modelling of the tensor $\sigma$, however, ad hoc assumptions and heuristics are employed. First of all, it is assumed that $\sigma$ can be written as
\begin{equation}\label{eqn_stresstens}
  \sigma = -\p I + \tau
\end{equation}
where $\p$ -- the pressure -- is a smooth scalar function and $\tau\colon \mathcal I \times \Omega \to \mathbb R^{3,3}$ is the tensor of shear stresses. 
It is assumed that $\tau(t,x)$ is symmetric and invariant under rigid body rotation. Furthermore, $\tau$  is a linear function of the velocity gradient $\nabla \v :=\bigl [\frac{\partial \v_j}{\partial x_i}\bigr]_{i,j=1,2,3}$.
\begin{remark}
  The separation of $\sigma$ into a pressure and shear stress part is motivated by requirement to recover the equations of hydrostatic or \emph{ideal flows}, where there are no shear stresses because of the absence either of motion or of viscosity. The symmetry of the tensor $\tau$ can be derived from the requirement that the angular momentum of a flow agglomerate is conserved. The invariance of $\tau$ with respect to rigid body rotations is derived from the general assumption that the considered fluid is isotropic, i.e., the physical properties of the fluid are the same in all spatial directions. Finally, the assumption that $\tau$ depends linearly on $\nabla u$ is well grounded for many fluids of interest. In fact, the validity of this assumption classifies a fluid as \emph{Newtonian fluid}.
\end{remark}
Under the given assumptions, $\tau$ is defined by the velocity field and two further parameters. It is commonly  written as 
\begin{equation}\label{eqn_stresstensor}
  \tau = \mu\bigl[\nabla \v + \nabla \v^\trp- \frac23(\dive \v)I\bigr] + \zeta (\dive \v)I.
\end{equation}
The parameter $\mu$ is the \emph{(first coefficient of the) viscosity} and has been experimentally determined and tabulated for many gases and fluids. The parameter $\zeta$ is called the \emph{bulk viscosity} and, in line with the \emph{Stokes Hypothesis}, often set to zero.

As for the left hand side in \eqref{eqn_mombalance}, one proceeds as follows. Let $\v_i$ denote the $i$-th component of $\v$, $i=1,2,3$. Then, an application of \eqref{eqn_RTT} gives that
\begin{equation}\label{eqn_RTTforvi}
  \ddt \int_{W} \rho \v_i \inva V = \int_W \pdt{\rho \v_i} + \dive (\rho \v_i \v) \inva V 
  = \int_W \pdt{\rho \v_i} + \rho \v_i \dive \v + \grad (\rho \v_i) \cdot \v \inva V 
\end{equation}
where basic vector calculus has been applied. If one considers $\nabla$ as the formal column vector of the three space derivatives, relation \eqref{eqn_RTTforvi} for all components of $v$ can be written in compact form as
\begin{equation}\label{eqn_RTTforv}
  \ddt \int_{W} \rho \v \inva V 
  = \int_W \pdt{\rho \v} + \rho \v \dive \v + (\v \cdot \nabla)(\rho \v)  \inva V .
\end{equation}
\subsubsection*{The Navier--Stokes equations}
In the preceding derivations, integral quantities over a flow agglomerate $W$ were considered. If the underlying \emph{continuity assumption} includes that all relations hold on arbitrary (small) control volumes $W$, instead of equating the integrals, one can equate the integrands pointwise in space, which leads to partial-differential equations.

Thus, putting together all assumptions and derivations, the balance of momentum \eqref{eqn_mombalance} for an isotropic Newtonian fluid under the \emph{Stokes Hypothesis} defines the (NSE) as
\begin{equation}\label{eqn_diffNSEgenform}
  \pdt{\rho \v} + \rho \v \dive \v + (\v \cdot \nabla)(\rho \v) 
  = \rho g - \nabla \p + \dive\bigl(\mu\bigl[\nabla \v + \nabla \v^\trp- \frac23(\dive \v)I\bigr]\bigr).
\end{equation} 
In the case that the flow is incompressible with a constant density $\rho \equiv \rho^*$ and a constant viscosity $\mu\equiv\mu^*$, the combination of \eqref{eqn_contieqint} and \eqref{eqn_diffNSEgenform} results in the system 
\begin{subequations}
\label{eqn_diffNSE}
\begin{align}
  \rho^* \bigl(\pdt{\v} + (\v \cdot \nabla) \v\bigr) + \nabla \p -\mu^* \Delta \v &= \rho^* g, \\
  \dive \v &= 0.
\end{align}
\end{subequations}
\begin{remark}
	Note that in the derivation, the pressure $\p$ is not a variable but a function which is assumed to be known and to describe the normal forces on a fluid element, cf.~\eqref{eqn_stresstens} and \cite[Ch.~1.1.ii]{ChoM93}. Thus, the NSE \eqref{eqn_diffNSEgenform} has only the velocity $\v$ as an unknown and the divergence free formulation \eqref{eqn_diffNSE} can be seen as an abstract ODE for $\v$ with an invariant \cite{HaiLW06}. Nonetheless, since it turns out that the divergence constraint defines the function $\p$, system \eqref{eqn_diffNSE} is commonly considered as an abstract differential-algebraic equation with $\v$ and $\p$ as unknowns.
\end{remark}
For the numerical treatment and for similarity considerations, one relates all dependent and independent variables to a characteristic length $L$ and characteristic velocity $V$ via
\begin{equation}\label{eqn_nondimvars}
	\v' = \frac \v V, \quad \p'=\frac{\p}{\rho^* U^2}, \quad 
	x' = \frac xL, \quad\text{and}\quad 
	t'  = \frac{tV}{L} \end{equation}
With this and with $f':=\frac{U^2}{L}g$, equation \eqref{eqn_diffNSE} can be rewritten in dimensionless form
\begin{subequations}\label{eqn_NSEnondim}
\begin{align}
	\frac{\partial }{\partial t'} \v' + (\v'\cdot \nabla)\v' - \frac 1\RE \Delta \v' + \grad \p' &= f', \\
	\dive \v' &= 0. \label{eqn_NSEnondim_divfree}
\end{align}
\end{subequations}
Thus, the system is completely parameterized by only one parameter $\RE:=\frac{UL\rho^*}{\mu^*}$,  the \emph{Reynolds number}. In what follows, we will always consider the dimensionless NSE \eqref{eqn_NSEnondim} but drop the dashes of the dimensionless variables.
\subsubsection*{Boundary conditions}
For the considered domain $\Omega$, let $\Gamma$ denote the boundary in the abstract and in the physical sense. If the domain is bounded by nonpermeable walls, then the \emph{no-slip} condition
\begin{equation}
	v = g \quad\text{on }\Gamma
\end{equation}
applies, where $g$ is the velocity of the wall. The no-slip conditions align well with experiments and macroscopic considerations; see~\cite[Ch.~5.3]{Lay08} for references but also for examples where no-slip conditions seem insufficient to describe the flow at walls.

If the domain is not fully bounded by walls, typically because the computational domain needs to be bounded whereas the physical domain of the flow is unbounded, \emph{artificial boundary} conditions are needed. We assume that the boundary is composed of a part $\Gammaw$ that is associated with a wall, a part $\Gammai$ where an (incoming) velocity profile is known, i.e., 
\begin{equation}
  \v = \v_i \quad\text{on }\Gammai,
\end{equation}
and a part $\Gammao$ that is often associated with an outflow.

There is no immediate physical insight into what should be the mathematical conditions at an artificial boundary condition $\Gammao$ that models the part where the fluid leaves the domain. Accordingly, outflow boundary conditions are motivated as being \emph{useful for the implementation of downstream boundary conditions} (see, e.g., \cite{Glo03, Pir89}) and equipped with the advice to put the outlet sufficiently far away from the region of interest; see~\cite[Ch.8.10.2]{FerP02}.

The most common outflow boundary conditions are the \emph{no stress} conditions:
\begin{equation*}
	\sigma n = 0 \quad\text{on }\Gammao,
\end{equation*}
where $\sigma$ is the stress tensor as defined in its general form in \eqref{eqn_stresstens}, the \emph{do-nothing} conditions:
\begin{equation}\label{eqn_bcs-donothing}
	\frac 1\RE \frac{\partial \v}{\partial n} - \p n = 0 \quad\text{on }\Gammao,
\end{equation}
that are formulated for the nondimensional Navier--Stokes equations \eqref{eqn_NSEnondim}, 
or the \emph{no gradient} conditions:
\begin{equation*}
	\frac{\partial \v}{\partial n} = 0 \quad\text{and}\quad\frac{\partial \p}{\partial n} = 0 \quad\text{on }\Gammao.
\end{equation*}
that are also applied at \emph{symmetry planes}.
\begin{remark}
The \emph{do-nothing} conditions have been extended; see~\cite{BraM14}, to the case of backflow, i.e., when some, possibly spurious, inflow occurs at the boundary $\Gammao$.
\end{remark}
\subsection{Formulation as operator DAE}\label{sect_NSE_opDAE}
In this subsection we provide yet another formulation of the dimensionless NSE \eqref{eqn_NSEnondim}, namely in the weak form as an operator DAE. This is a DAE in an abstract setting, where the solution is an element of a Sobolev space instead of a vector. 
We consider here homogeneous Dirichlet boundary conditions. More general boundary conditions with in- and outflow may modeled in a similar way. 
We define the spaces 
\[ 
  \V := [H^1_0(\Omega)]^n,\quad  
  \cH:=[L^2(\Omega)]^n,\quad \text{and } 
  \Q:= L^2(\Omega)/\R. 
\]
By $\V'$ we denote the dual space of $\V$. 
Note that the spaces $\V$, $\cH$, $\V'$ form a Gelfand or evolution triple \cite[Ch.~23.4]{Zei90a}.  
Furthermore, we define $\calW(0,T)$ as the space of functions $u \in L^2(0,T;\V)$, which contain a weak time derivative $\dot u \in L^2(0,T;\V')$. 
Well-known embedding results then imply $u \in C([0,T];\cH)$, cf.~\cite[Lem.~7.3]{Rou05}.

We now consider the weak formulation of \eqref{eqn_NSEnondim} in operator form. 
This means that for given right-hand sides $\calF \in L^2(0,T; \V')$, $\calG \in L^2(0,T;\Q')$ and an initial condition $\inicond \in \cH$, we seek for a pair $(\v, \p)\in \calW(0,T)\times L^2(0,T;\Q)$ satisfying 
\begin{subequations} \label{eqn_nseBlankWeak}
\begin{align}
  \dot \v(t) + \mathcal K(\v(t)) - \mathcal B'\p(t) &= \calF(t) \phantom{\calG\inicond} \text{ in } \V',\\
  \mathcal B \v(t) \phantom{i + \calB'\p(t)} &= \calG(t) \phantom{\calF\inicond} \text{ in } \Q', \label{eqn_nseBlankWeak_b} \\
  \v(0) &= \inicond \phantom{\calG\calF(t)} \text{ in } \cH
\end{align}
\end{subequations}
a.e. on $(0,T)$. The derivative of $\v$ should be understood in the weak sense. 
The operators $\mathcal K\colon \V \to \V'$ and $\mathcal B\colon \V \to \Q'$ are defined via
\begin{equation} \label{eqn_defK}
  \langle \mathcal K(\v),w \rangle 
  = \int_\Omega (\v\cdot \nabla)\v \cdot w\dx + \frac{1}{\RE} \int_\Omega \nabla\v \cdot \nabla w\dx
\end{equation}
and 
\begin{equation}
  \langle \mathcal B\v, q \rangle  
  = \int_\Omega (\ddiv \v)q \dx
  = \langle \v, \calB' q \rangle,
\label{eqn_defB}
\end{equation}
respectively, for a given $\v\in \V$ and for all test functions $w \in \V$ and $q \in \Q$. 
We emphasize that system \eqref{eqn_nseBlankWeak} not only covers the NSE but also more general flow equations, since we have introduced an inhomogeneity $\calG$.  

For results on the existence solutions to the weak formulation of the NSE, we refer to \cite{Tar06}. For a compact summary that also considers nonzero $\mathcal  G$ in \eqref{eqn_nseBlankWeak_b} see~\cite{AltH15}. The differential-algebraic structure of \eqref{eqn_nseBlankWeak} and the possible decoupling of differential and algebraic parts and variables has been discussed in~\cite{Hei14}.
 \section{Semi-discrete equations}\label{sect_dae} In this section, we consider the DAE, which results from a spatial discretization of system \eqref{eqn_NSEnondim} or \eqref{eqn_nseBlankWeak}. 
Using finite elements, finite differences, or finite volumes, we obtain a DAE of the form 
\begin{subequations}
\label{eqn_NS_disc}
\begin{align}
  M \dot v + K(v) - B^T p &= f, \label{eqn_NS_disc_a}\\
  B v &= g. \label{eqn_NS_disc_b}
\end{align}
\end{subequations}
Here, $v$ and $p$ denote the finite-dimensional approximations of $\v$ and $\p$, respectively. 
In addition, we assume a consistent initial condition of the form $v(0)=\qo$, i.e., we demand $B\qo = g(0)$. 

There are many reasons to include a right-hand side $g$ in equation \eqref{eqn_NS_disc_b} rather than a zero as it seems naturally for incompressible flows. Firstly, a nonzero $g$ in the continuous continuity equations \eqref{eqn_NSEnondim_divfree} may appear also in generalizations of the NSE model to fluid-structure interactions (see, e.g., \cite{Ray10}) or in optimal control setups; see, e.g., \cite{Hin00}. Secondly, a nonzero $g$ maybe a numerical artifact from the semidiscretization of flows with nonzero Dirichlet boundary conditions as in the example we will provide below. Finally, in view of analysing the DAE, this inclusion of a nonzero $g$ leads to a better understanding as one can track where the derivatives of the right-hand sides appear within the solution. This is of importance also in the case $g=0$, since inexact solves lead to errors in the \eqref{eqn_NS_disc_b} that act like a nonsmooth inhomogeneity.
\subsection{Spatial discretization by finite elements}\label{sect_dae_discretization}
We consider a discretization by finite elements, i.e., we construct finite-dimensional (sub)spaces $V_h$ and $Q_h$ of $\V$ and $\Q$, respectively, based on a (shape) regular triangulation $\calT$ of the polygonal Lipschitz domain $\Omega$, cf.~\cite{Cia78}. 
Given a basis $\{\varphi_1, \dots, \varphi_n  \}$ of $V_h$, we can identify the finite-dimensional approximation of the velocity $\v(t)$ by the coefficient vector $v(t)\in\R^n$. 
The discrete representative of the pressure $\p$(t) is denoted by $p(t)\in\R^m$ and corresponds to a basis $\{ \psi_1, \dots, \psi_m \}$ of $Q_h$.  

With the basis functions of $V_h$ we define the symmetric and positive definite mass matrix $M \in \R^{n,n}$ by 
\[
  M := [m_{jk}] \in \mathbb R^{n\times n},\qquad 
  m_{jk} := \int_\Omega \varphi_j \cdot \varphi_k \dx.
\]
For specific discretization schemes the mass matrix may even equal the identity matrix. 
The discretization of the nonlinearity of the Navier-Stokes equation, i.e., the discretization of the operator $\calK$ in Section~\ref{sect_NSE_opDAE}, is denoted by $K\colon \R^n \to \R^n$. 
Note that the given model also includes linearizations of the Navier-Stokes equation such as the unsteady Stokes or Oseen equation. In this case, $K$ can be written as a $n\times n$ matrix.  
In view of the index analysis of system \eqref{eqn_NS_disc}, however, this is not of importance; see~Appendix~\ref{app_sindex}. 
We define for $v\in \R^n$ and its representative $\tilde v = \sum_{j=1}^n v_j \varphi_j \in V_h$, 
\begin{equation}\label{eqn_defdiscreteK}
  K_j(v) := \int_\Omega (\tilde v \cdot \nabla) \tilde v \cdot \varphi_j \dx 
     + \frac{1}{\RE} \int_\Omega \nabla \tilde v \cdot \nabla \varphi_j \dx. 
\end{equation}
Finally, we define the matrix $B$, which corresponds to the divergence operator, i.e.,   
\[
  B := [b_{ij}]\in\R^{m\times n},\qquad
  b_{ij} = \int_\Omega \psi_i \ddiv \varphi_j \dx.
\]
Since $\ddiv$ and $-\nabla$ are dual operators in the continuous setting, in the semi-discrete equations we get a matrix $B$ and (minus) its transpose $B^T$. 
Note that this yields the saddle point structure of system \eqref{eqn_NS_disc}.

Especially for the stable approximation of the pressure, it is necessary that the chosen finite element spaces are compatible \cite[Ch.~VI.3]{BreF91}. 
Let $V_h$ and $Q_h$ denote again finite-dimensional (sub)spaces of $\V$ and $\Q$, respectively, with 
\[
  \dim V_h = n, \qquad
  \dim Q_h = m < n.
\]
The spaces $V_h$ and $Q_h$ are compatible if they satisfy a so-called {\em inf-sup} or {\em Ladyzhenskaya-Babu${\check {\text s}}$ka-Brezzi} condition \cite[Ch.~VI.3]{BreF91}. 
This means that there exists a constant $\beta>0$, independent of the chosen mesh size, such that
\[
 \adjustlimits \inf_{p_h\in Q_h} \sup_{v_h \in V_h} 
 \frac{| \langle \ddiv v_h, p_h\rangle |}{\Vert v_h\Vert_\V \Vert p_h\Vert_\Q}
 \geq \beta.
\]
The inf-sup condition, with $\beta$ independent of $h$, is a necessary condition for the convergence of the FEM; see, e.g., \cite{Lad69} or \cite{GirR86}. For a fixed spatial discretization, this condition implies that the matrix $B$ resulting from the discretization scheme $V_h$, $Q_h$ is of full rank. One scheme which is known {\em not} to satisfy this condition (in dimension $d=2$) is given by 
\[
  V_h = [\calP_1(\calT) \cap \V]^2, \qquad
  Q_h = \calP_0(\calT)/\R. 
\]
Here $\calP_k(\calT)$ denotes the space of piecewise polynomials of degree $k$. Thus, the space $V_h$ involves continuous finite element functions of polynomial degree one, whereas $Q_h$ consists of piecewise constant functions. Furthermore, $V_h$ includes homogeneous Dirichlet boundary conditions, since they are part of $\V$. It has been shown that that the velocity space is too small and thus $\beta \not > 0$ \cite[Ex.~VI.3.1]{BreF91}.  
This flaw can be fixed in various ways. One possibility is given by introducing so-called {\em nonconforming} finite element spaces. Crouzeix and Raviart \cite{CroR73} proposed to use piecewise polynomials of degree one, which are only continuous in the midpoints of the edges of the triangulation. Together with homogeneous Dirichlet boundary conditions the ansatz space is denoted by $\CR_0(\calT)$. 

Another strategy is to enrich the space $V_h$ by so-called {\em bubble functions}, namely $\calB$, cf.~\cite{BerR85}. With this we introduce ansatz functions of polynomial degree two. 
A small collection of stable finite element schemes is given in Table~\ref{tab_stableFE}. Further stable schemes are addressed in \cite[Ch.~II]{GirR86} as well as in \cite[Ch.~3]{GreS00}.

\begin{remark}
Factoring out $\mathbb R$ in the definition of the pressure approximation spaces $\mathcal Q_h$ is needed to eliminate the kernel of the gradient operator or $B^T$, which are the constant functions both in finite and infinite dimensions; see \cite{GirR86}. This reflects the fact that in internal flows the pressure is defined up to a constant.  
In fact, the standard results for inf-sup stability always assume the case of internal flow, i.e. zero Dirichlet conditions everywhere at the boundary. There have been but a few attempts to derive inf-sup conditions that apply to nonzero boundary conditions \cite{GunH92}.
\end{remark}

\begin{table}
\label{tab_stableFE}
\caption{Selective overview of stable finite element schemes for flow equations.}
\begin{tabular}{c|l|l} 
  discretization scheme & discrete function spaces      & inf-sup stability \\ \hline\hline
Taylor-Hood             & $V_h = [\calP_k(\calT) \cap \V]^d,\ \ k\ge 2 $      & \cite{TayH73,Ver84b} \\ 
(continuous pressure)   & $Q_h = [\calP_{k-1}(\calT) \cap \V] /\R$   & \\ \hline
Crouzeix-Raviart        & $V_h = [\CR_0(\calT)]^d$      & \cite{CroR73,BecM11}\\ 
(discontinuous velocity)& $Q_h = \calP_0(\calT)/\R$     &  \\ \hline
Bernardi-Raugel         & $V_h = [\calP_1(\calT) \cap \V]^d \oplus \calB$      & \cite{BerR85} \\ 
(edge-bubble functions) & $Q_h = \calP_0(\calT)/\R$     & \\ \hline
Rannacher-Turek      & on quadrilateral meshes              & \cite{RanT92} 
\end{tabular}
\end{table}
\subsection{Finite volumes and finite differences}\label{sec_FVMFDM}
In this section we briefly touch the spatial discretization of the incompressible NSE by the methods of \emph{finite volumes} (FVM) and \emph{finite differences} (FDM). If applied to \eqref{eqn_NSEnondim} in a straight forward manner, both approaches lead to a DAE of type \eqref{eqn_NS_disc}. 

Finite volume approximations base on the integral formulation of the conservation laws as it reads for the momentum equation
\begin{equation}\label{eqn_intforFVM}
	\int_W \pdt{\rho \v} \inva V + \int_{\partial W} \rho \v \v \cdot n\inva S 
	=- \int_{\partial W} \p n \inva S + \int_{\partial W} \tau \cdot n \inva S + \int_W \rho g \inva V,
\end{equation}
cf.~\eqref{eqn_mombalance} and \eqref{eqn_stresstens}. Balances, like \eqref{eqn_intforFVM} in particular, hold for the whole domain of computation. For the discretization, the domain of the flow is subdivided into small volumes often referred to as cells. The velocities and the pressure are assumed to be, say, constant over the cells, and their approximated values are determined by evaluating and equating the volume and surface integrals associated with every cell; see~\cite[Ch. 8.6]{FerP02}.

Because of its flexibility in the discretization, because of its variants that provide unconditional stability properties, and since typical models for the effective treatment of turbulence are formulated as conservation laws too, the FVM is the method of choice in most general purpose solvers. As for the discussion in the DAE context we note that, typically, the DAE \eqref{eqn_NS_disc} is never assembled but rather decoupled during the time discretization; see~\cite[Ch.~7]{FerP02} and Section~\ref{sec_SIMPLE} below.

Approximations of the NSE via FDM are not widely used because of the known flaws of finite difference approximations like high regularity requirements and confinement to regular grids. Nonetheless, FDM discretizations have been successfully used in flow discretizations and are probably still in use in certain specified, say, single-purpose codes. Also, FDM are often well suited for teaching the fundamentals of flow simulation of laminar and turbulent setups; see~\cite{GriDN97}.
\subsection{Index of the DAE}\label{sect_dae_rankdeficit} 
As argued above, in \eqref{eqn_NS_disc} the coefficient matrix $B$ may be rank-deficient and, thus, make the \emph{differentiation index} not well-defined \cite{Wei96ppt}. 
However, the \emph{strangeness index}~\cite{KunM94, KunM01, KunM06} that applies to over- and underdetermined DAE systems and, thus, also can be determined in the case of a rank-deficient coefficient matrix $B$.
A rough index analysis has been realized in \cite{Wei97} with the result that under general and reasonable assumptions, system \eqref{eqn_NS_disc} has strangeness index 1. See also Appendix~\ref{app_sindex}, where the strangeness index has been determined in a rigorous way.
This matches the, say, observations in \cite[Ch.~VII.1]{HaiW96} that the system is of (differentiation) index 2 if the matrix $B$ is of full rank.

In the sequel we always assume that the discretization scheme is chosen in such a way that the matrix $B$ is of full rank. Note that this can always be realized, choosing linearly independent finite element basis functions with respect to the space $\Q$, i.e., keeping in mind that constant functions are in the same equivalence class as the zero-function. 
If the ansatz functions form a partition of unity, then this means nothing else than eliminating one of these ansatz functions.  
Numerical schemes used in practice usually satisfy this condition as the discrete ansatz space for the pressure $Q_h$ is chosen appropriately. 
 \section{Fully discrete approximation schemes}\label{sect_solution}
The final step of a numerical approximation of the infinite-dimensional NSE \eqref{eqn_NSEnondim} is the numerical time integration of the semi-discrete approximation \eqref{eqn_NS_disc}. For this, one discretizes the time interval $(0, T)$ via the grid 
\begin{equation*}
  t_0 :=0 < t_1 < t_2 < \dotsm < t_N = T
\end{equation*}
and computes a sequence $(v^k , p^k)_{k\in \mathbb N}$ of values that are supposed to approximate the solution of \eqref{eqn_NS_disc} at the discrete time instances, i.e., $v^k \approx v(t^k )$ and $p^k  \approx p(t^k )$.

For example, a single-step time integration scheme applied to an ODE $\dot x = \phi(t, x)$ leads to an approximating sequence defined through
\begin{equation}\label{eq_fdiscODE}
	x^{k+1}  = \Phi (t^k , x^k , \tau^k ),
\end{equation}
where $\tau^k := t^{k+1} - t^k $ and $\Phi$ is the increment function of the single-step scheme.
\subsection{Index and causality of discrete systems}\label{sect_solution_index}
For the analysis, we will restrict our considerations to a linear time-invariant setup. This is no restriction, since in the particular semi-linear semi-explicit form, the DAE structure is not affected by the nonlinearity the more that in the presented methods the momentum equation is typically discretized with an explicit scheme. Also we assume that the time grid is equidistantly spaced and of size $\tau=t^{k+1} - t^k $.

We will cast the fully-discrete schemes into the standard form
\begin{equation}\label{eqn_EAstandardform}
  \mathcal Ex^{k+1} = \mathcal A x^k + f^k,
\end{equation}
with coefficient matrices $\mathcal E$ and $\mathcal A$ and $x^k$ containing all variables like $x^k = [v^k ; p^k]$. 
We call such a fully discrete approximation scheme \eqref{eqn_EAstandardform} a {\em difference-algebraic equation} ($\Delta$AE).
\begin{remark}
This recast of the scheme, which will basically amount to a shift of the time indices for some variables, is needed to avoid ambiguities. In fact, for continuous time DAEs with delays, one can deliberately apply the operations of time-shift and differentiation to produce equivalent formulations of different indices \cite{Ha15}. 
\end{remark}
In what follows, we argue that the matrix pair \EApair~can be assumed to be regular (cf.~\cite[Def. 2.5]{KunM06}) such that it can be brought into a particular Kronecker form 
\begin{equation}\label{eq:eakronform}
	\EApair \sim 
	\bigl (
	\begin{bmatrix}
		I & 0 \\ 0 & N
	\end{bmatrix},
	\begin{bmatrix}
		J & 0 \\ 0 & I
	\end{bmatrix}
	\bigr),
\end{equation}
where $J$ is a matrix in Jordan form and $N$ is a nilpotent matrix, cf.~\cite[Thm.~2.7]{KunM06}. 
\begin{definition}\label{def_disc_kron_index}
The {index of nilpotency} -- that integer $\nu$ for which $N^\nu=0$ while $N^{\nu-1}\neq 0$ -- is called the \emph{index of the matrix pair} \EApair; see~\cite[Def.~2.9]{KunM06}, or the \emph{Kronecker index} of the $\Delta$AE \eqref{eqn_EAstandardform}, cf.~\cite[p.~39]{Rei06}.
\end{definition}
The excursion to the discrete Kronecker index is necessary, since the common time stepping schemes considered in Section \ref{sect_solution_methods} below can not be interpreted as a time-discretization of time-continuous DAE. Equivalently, we can can classify the discrete schemes using the notion of \emph{causality} that attributes systems with states depending only on the past or current inputs.
\begin{definition}[{\cite[Def.~8-1.1]{Dai89}}]
The sequence $x^k $, $k=1,2,\dots$, defined through \eqref{eqn_EAstandardform} is called \emph{causal}, if $x^k$ is determined completely by an initial condition $x^0$ and former (and the current) inputs $f^0, f^1, \dotsc, f^k$.
\end{definition}

Conversely, in a \emph{noncausal} system, the current state $x^k$ depends on future inputs like $f^{k+1} $. We will analyse time-stepping schemes for causality, which in the considered case is equivalent to being of Kronecker index $1$ \cite[Thm.~8-1.1]{Dai89} and discuss how and why a noncausal system poses difficulties in the numerical approximation. 

To introduce the procedure, to fix the notation, and to have a benchmark for further comparisons, we start with analysing an \emph{implicit-explicit Euler} discretization of a linearized version of \eqref{eqn_NS_disc}, namely
\begin{subequations} \label{eqn_diffschmEuler}
\begin{align}
	\tfrac 1\tau Mv^{k+1} &= (\tfrac  1\tau M+A)v^k + B^T p^k + f^k, \label{eqn_diffschmEuler_mom}\\
	Bv^{k+1} &= g^{k+1}. \label{eqn_diffschmEuler_conti}
\end{align}
\end{subequations}
This defines the difference equations for $v^k$ and $p^k$, which approximate the velocity and pressure at the discrete time instances $t^k$, $k=1,2,\dots$. Similarly, $f^{k}$ and $g^{k}$ stand for the approximations of $f$ and $g$ at the time instances.

The difference scheme is of the standard form \eqref{eqn_EAstandardform} with $x^k = [v^k ; p^k ]$, $h^k = [f^k ; g^k ]$, a shift of the index in \eqref{eqn_diffschmEuler_conti}, and 
\begin{equation}\label{eqn_EAplainEul}
	\bigl(\mathcal E, \mathcal A\bigr) = 
	\bigl( \begin{bmatrix} \tfrac 1\tau M & 0 \\ 0 & 0 \end{bmatrix},  
	       \begin{bmatrix} \tfrac 1\tau M+A & B^T \\ B & 0 \end{bmatrix} \bigr).
\end{equation}
In Appendix~\ref{app_diffAE} it is shown that under standard conditions and for $\tau$ sufficiently small, the pair $(\mathcal E, \mathcal A)$ in \eqref{eqn_EAplainEul} is regular and equivalent to 
\begin{equation*}
	\bigl( \begin{bmatrix} I & 0& 0 \\ 0 & 0& 0\\ 0 & I& 0 \end{bmatrix},  
	       \begin{bmatrix} * & 0& 0 \\ 0& I & 0 \\0& 0 & I \end{bmatrix} \bigr).
\end{equation*}
This means that the difference scheme based on an \emph{implicit-explicit Euler} discretization is of \emph{Kronecker index} 2. In fact, a straight forward calculation reveals that in \eqref{eqn_diffschmEuler} the state $p^k$ depends on $g^{k+1}$.  
\begin{remark}
The index $k$ for the pressure $p$ in \eqref{eqn_diffschmEuler} is consistent with $f^k$, as can be directly derived from the case that $A=0$ and $g=0$. A fully implicit scheme, i.e., considering $f^{k+1} $ and $p^{k+1}$ in \eqref{eqn_diffschmEuler_mom}, can not be brought into the standard form \eqref{eqn_EAstandardform}.
\end{remark}

\subsection{Inherent instabilities of \dDAE s of higher index}
In this section, we illustrate a mechanism that leads to a numerical instability and, thus, possibly to divergence of the approximation of a dynamical system through a time discretization, i.e. a \dDAE, with an index greater than $1$. Consider a DAE in Kronecker form 
\begin{equation*}
	N\dot x = x + g
\end{equation*}
with $N\neq0$ and $N^2=0$ and it's time discrete approximation through an Euler scheme:
\begin{equation}\label{eqn_ind2ddae}
	\frac 1\tau Nx^{k+1}  = (\frac 1\tau N + I)x^k  + g^k .
\end{equation}
The \dDAE~\eqref{eqn_ind2ddae} is of Kronecker index 2, according to Definition \ref{def_disc_kron_index} and as it can be read off after a premultiplication by $(I-N)$, and it has the solution
$$
x^k = -g^k - \frac 1\tau(Ng^{k+1} - Ng^k ),
$$
as it follows from an adaption of the arguments in \cite[Lem. 2.8]{KunM06} to the discrete case. From this solution representation one can conclude, that the solution to a \dDAE\ of higher index that discretizes a DAE may depend on numerical differentiations and that any error in the computation may be amplified by the factor $\tau^{-1}$. Note that for index-1 \dDAE s, where $N=0$, this derivative is not present and that for even higher indices higher (numerical) derivatives will appear in the solution. Also note, that for systems that are not in Kronecker form such as \eqref{eqn_diffschmEuler}, these derivations will be realized implicitly; see~\cite{AltH15}.
\subsection{Common time-stepping schemes as index-1 $\Delta$AEs}\label{sect_solution_methods} In this subsection, we discuss the different strategies, which are used in practice, to solve the spatially discretized NSE \eqref{eqn_NS_disc}. 
In practical applications one typically uses schemes that decouple pressure and velocity computations. Although, as we have argued from a DAE perspective \cite{AltH15}, this may lead to instabilities, the advantages that
\begin{itemize}
	\item one has to solve two smaller systems rather than one large and
	\item one basically solves Poisson equations and convection-diffusion rather than saddle-point problems \end{itemize}
seem to prevail. The second point is important in so far as for the automated solution of the coupled problem there are no generally well performing solvers for the arising linear systems. For example, in the cases of practical relevance namely that of high Reynolds numbers or highly non-normal coefficient matrices, there are no efficient preconditioners for solving the resulting saddle-point systems iteratively. 
In fact, recent developments in this direction were tested in low to moderate $\RE$-setups and the reported numbers indicate an unfavourable scaling of the iteration numbers with \RE; see, e.g., \cite[Ch.~8]{ElmSW05} or \cite{HeV16}. Sometimes, in particular when it comes to parallelization, \emph{Gauss--Seidel}-type methods are used, which iteratively and locally update the variables. As far as convection and diffusion is concerned such a formulation, that the quantities of interests in a cell depend on the local neighbors, is physically well backed. 

Another common feature of the schemes used in practice is their explicit approach to the momentum equation which avoids the repeated assembling of Jacobians. The error is then either controlled through a small time-step or through some fixed point iterations.

In this subsection, we consider the schemes \emph{Projection} as it was described in \cite{Gre90}, \emph{SIMPLE} as in 
\cite{FerP02}, and \emph{artificial compressibility} \cite{FerP02}. We will use a description and formulation general enough, to also accommodate numerous variants of the methods.

\subsubsection{Projection}
The principle of these methods is to solve for an intermediate velocity approximation that does not need to be divergence-free, and project it onto the divergence-free constraint in a second step. As far as the velocity approximation is concerned, \emph{Projection} methods can be formulated both in infinite and finite dimensions. Since the first work by Chorin \cite{Cho68}, a number of variants have been developed mainly proposing different approaches to the approximation of the pressure tackling or circumventing the need of solving a \emph{Poisson equation} for the pressure, which requires certain regularity assumptions (cf.~\cite[p.~642]{GreS00}). Another problem is the requirement of boundary conditions for the pressure update that do not have a physical motivation and may cause inaccuracies close to the boundary. Nonetheless, these schemes have been extensively studied and certain heuristics ensure satisfactory convergence behaviour; see, e.g, \cite{Gre90,GreS00}. 

As an example for a \emph{Projection} scheme, we present the variant proposed in \cite{Gre90}:
\begin{enumerate}
	\item Solve for intermediate velocity with the old pressure
        \begin{equation}\label{eqn_tildev}
   	\tfrac  1\tau M\tilde v^{k+1} = (\tfrac 1\tau M+A)v^k + B^T p^k + f^k.
        \end{equation}
\item Determine the new velocity $v^{k+1}$ as the projection of the intermediate velocity onto $\ker B$ by solving  
	\begin{subequations}\label{eqn_diffschmPrj_prjstep}
		\begin{align}
			Mv^{k+1} - B^T\phi^{k+1} &=M\tilde v^{k+1}, \\
			Bv^{k+1} &= g^{k+1}. \label{eqn_diffschmPrj_conti}
		\end{align}
	\end{subequations}
	\item Update the pressure via
		\begin{equation}\label{eqn_prjctn2pupd}
		p^{k+1} = p^k + \frac 2\tau \phi^{k+1}.
	\end{equation}
\end{enumerate}
\begin{remark}
Instead of solving the saddle-point problem \eqref{eqn_diffschmPrj_prjstep} as a whole and in order to avoid the division of a numerically computed quantity by $\tau$ in \eqref{eqn_prjctn2pupd}, one can decouple the system. Then, one solves for $\tilde \phi^{k+1} := \frac 2\tau \phi^{k+1} $  through 
	$$- BM^{-1}B^T\phi^{k+1} =\frac 2\tau (B\tilde v^{k+1} - g^{k+1}) $$
	and obtains the updates via
	$$v^{k+1} = \tilde v^{k+1} + \frac \tau 2 M^{-1}B^T\tilde \phi ^{k+1}$$
	and 
	$$ p^{k+1} =p^k + \tilde \phi^{k+1}. $$
\end{remark}
\begin{remark}
The formula for the update of the pressure is derived from the relation 
\[ \phi(t+\tau) = -\frac{\tau^2}{2}\dot p(t) + \mathcal O(\tau^3), \]
cf.~\cite[p.~595]{Gre90}, that holds under certain regularity assumptions. Note also that, if formulated for the space-continuous problem, this \emph{Projection 2} algorithm requires a sophisticated treatment of the boundary conditions. The presented variant is a simplification for spatially discretized equations; see~\cite[Ch.~3.16.6c]{GreS00}.
\end{remark}

As a single system, this projection scheme defines a \dDAE~in  the form of \eqref{eqn_EAstandardform} with $x^k := [\tilde v^k ; \phi^k ; v^k ;p^k ]$ and
\begin{equation}\label{eqn_eapair_prjctn}
	\bigl(\mathcal E, \mathcal A\bigr) = 
	\bigl( 
	\begin{bmatrix} \tfrac 1\tau M & 0 & 0 & 0\\
									0 & 0 & 0 & 0 \\
									-M & -\frac \tau 2 B^T & M & 0 \\
									0 & -I & 0 & I

	\end{bmatrix},  
	\begin{bmatrix} 0 & \tfrac 1\tau M+A & 0 & B^T \\
									-\frac 2\tau B & -BM^{-1}B^T & 0 & 0 \\
									0 & 0 & 0 & 0 \\
									0 & 0 & 0 & I
	\end{bmatrix}
	\bigr),
\end{equation}
which is a matrix pair of Kronecker index 1.
\subsubsection{SIMPLE scheme -- implicit pressure correction}\label{sec_SIMPLE}
The SIMPLE scheme and its variants are based on the decomposition 
\begin{equation}\label{eqn_simplecorrs}
v^{k+1} = \tilde v^{k+1} + v_\Delta ^{k+1} \quad\text{and}\quad p^{k+1} = p^k + p_\Delta^{k+1} ,
\end{equation}
where $\tilde v^{k+1}$ is the tentative velocity computed by means of the old pressure. We present the basic variant, in which the velocity correction $v_\Delta^{k+1} $ is discarded when solving for $p_\Delta^{k+1} $:
\begin{enumerate}
	\item Solve for the intermediate velocity $\tilde v^{k+1}$ with the old pressure as in \eqref{eqn_tildev}. 
	\item Compute $p_\Delta^{k+1} $ through
	\begin{equation}
	B(\tfrac 1\tau M+A)^{-1}B^Tp_\Delta^{k+1} =  	-B\tilde v^{k+1} + g^{k+1} \label{eqn_diffschmSMPL_pcpe}
     \end{equation}
as the correction to $p^k$ such that the
	\item updates of the velocity and pressure defined through
\begin{subequations} \label{eqn_diffschmSMPL}
\begin{align}
	v^{k+1} &= \tilde v^{k+1}  + (\tfrac 1\tau M+A)^{-1}B^Tp_\Delta^{k+1}, \label{eqn_diffschmSMPL_velpc}\\
	p^{k+1} &= p^k + p_\Delta^{k+1}. 
\end{align}
\end{subequations}
jointly fulfill the time discrete momentum and the continuity equation. 
\end{enumerate}
\begin{remark}
Step (2) of the presented SIMPLE algorithm computes $p_\Delta^{k+1} $ under the temporary assumption that $v_\Delta^{k+1} =0$. This is hardly justified and probably a reason for slow convergence, when applied in a fixed-point iteration within fully implicit schemes. Certain variants of the SIMPLE scheme try to approximate $v_\Delta^{k+1} $ at this step, e.g., through interpolation; see~\cite[Ch.~7.3.4]{FerP02}.

\end{remark}

As a single system, the SIMPLE scheme defines a \dDAE~in  the form of \eqref{eqn_EAstandardform} with $x^k := [\tilde v^k ; p_\Delta^k ; v^k ;p^k ]$ and
\begin{equation}\label{eqn_eapair_SMPL}
	\bigl(\mathcal E, \mathcal A\bigr) = 
	\bigl( 
	\begin{bmatrix} \tfrac 1\tau M & 0 & 0 & 0\\
									0 & 0 & 0 & 0 \\
									I & (\tfrac 1\tau M+A)^{-1}B^T & -I & 0 \\
									0 & -I & 0 & I

	\end{bmatrix},  
	\begin{bmatrix} 0 & \tfrac 1\tau M+A & 0 & B^T \\
								- B & -B(\frac 1\tau M + A)^{-1} B^T & 0 & 0 \\
									0 & 0 & 0 & 0 \\
									0 & 0 & 0 & I
	\end{bmatrix}
	\bigr),
\end{equation}
which is a matrix pair of Kronecker index 1; see~Appendix \ref{sec_appdx_kronind_SMPL}.
\subsubsection{Artificial compressibility} \label{sec_artificialComp}
In this class of methods, the divergence-free constraint \eqref{eqn_NSEnondim_divfree} is relaxed by adding a scaled time-derivative of the pressure, i.e., 
\begin{align}\label{eqn_artcompdivfree}
	\frac 1\beta\pdt \p +  \dive \v &= 0.
\end{align}
The parameter is assumed to satisfy $\beta \gg 1$. The corresponding spatial semi-discretization \eqref{eqn_NS_disc} then reads
\begin{subequations}
\label{eqn_artcompsemdisc}
\begin{align}
  M \dot v + K(v) - B^T p &= f, \\
  M_p\dot p + B v &= g  \label{eqn_artcompsemdisc_compconstr} 
\end{align}
\end{subequations}
with $M_p$ denoting the mass matrix of the pressure approximation. 

In theory, system \eqref{eqn_artcompsemdisc} is an ODE and could be solved by standard time-stepping schemes. In practice, however, the common solution approaches to System \eqref{eqn_artcompsemdisc} decouple pressure and velocity computations through introducing auxiliary quantities and, thus, a DAE structure. A possible method, that defines a pressure update similar to the SIMPLE scheme \eqref{eqn_diffschmSMPL}, can be interpreted and implemented as follows, cf.~\cite[Ch.~7.4.3]{FerP02}. As in the SIMPLE approach a first velocity approximation $\tilde v^{k+1} $ is computed via \eqref{eqn_tildev} on the base of the old pressure value $p^k$. 

Next, one substracts the contribution of the old pressure, 
\begin{equation}\label{eqn_corrthev}
	\bar v^{k+1} := \tilde v^{k+1} - (\frac 1\tau M + A)^{-1}B^Tp^k
\end{equation}
and defines the new velocity $v^{k+1}$ through the linear expansion of $\bar v^{k+1} $ in terms of the pressure gradient, i.e., 
\begin{equation}\label{eqn_vkk_pc_ansatz}
	v^{k+1} = \bar v^{k+1} + \frac{\partial (\bar v^{k+1} )}{\partial (B^Tp)}[B^Tp^{k+1} - B^Tp^k ] = \bar v^{k+1} + (\frac 1\tau M + A)^{-1}B^Tp_\Delta^{k+1}.
\end{equation}
Note that one uses \eqref{eqn_corrthev} to determine the Jacobian $\frac{\partial (\bar v^{k+1} )}{\partial (B^Tp)}$ and that we have defined $p_\Delta^{k+1} = p^{k+1} -p^k $. 
Expression \eqref{eqn_vkk_pc_ansatz} is then inserted in a, e.g., \emph{implicit Euler} discretization of \eqref{eqn_artcompsemdisc_compconstr}, which gives an equation for $p_\Delta ^{k+1} $. Accordingly, the fully discrete scheme using artificial compressibility reads 
\begin{subequations} \label{eqn_diffschm_artcomp}
\begin{align}
	(\tfrac 1\tau M+A)\tilde v^{k+1} &= \tfrac  1\tau Mv^k + B^T p^k + f^k, \label{eqn_diffschm_artcomp_tv}\\
	\bar v^{k+1} &=\tilde v^{k+1} - (\tfrac 1\tau M+A)^{-1}B^Tp^k  , \\
	\frac 1{\beta \tau} p_\Delta^{k+1} + B\bar v^{k+1} + B(\tfrac 1\tau M+A)^{-1}B^Tp_\Delta^{k+1} &= g^{k+1} , \label{eqn_diffschm_artcomp_pcpe}\\
	v^{k+1} &= \bar v^{k+1}  - (\tfrac 1\tau M+A)^{-1}B^Tp_\Delta^{k+1}, \label{eqn_diffschm_artcomp_velpc}\\
	p^{k+1} &= p^k + p_\Delta^{k+1}.
\end{align}
\end{subequations}
The corresponding \dDAE equals the \dDAE of the SIMPLE scheme up to a slight modification of the equation for the pressure update such that we can use the arguments laid out in Appendix \ref{sect_solution_indexred} to conclude that it is of index 1.
\begin{remark}
Note that in practice there may be nonlinear solves to determine, e.g., $\tilde v^{k+1}$ such that, e.g., the correction $\bar v^{k+1} $ is possibly different from $\tilde v^{k+1}$ computed with $p^k =0$.
\end{remark}

\subsection{Time-stepping schemes resulting from index reduction}\label{sect_solution_indexred}
Besides the presented schemes in Section~\ref{sect_solution_methods}, one may also apply an index reduction to system~\eqref{eqn_NS_disc} and then discretize in time. This then also leads to matrix pairs $(\calE, \mathcal A)$ of Kronecker index 1.  
\subsubsection{Penalty methods}\label{sect_solution_indexred_penalty}
Similar to Section~\ref{sec_artificialComp}, we may reduce the index of the DAE~\eqref{eqn_NS_disc} by relaxing the divergence-free constraint or, in other words, add a penalty term~\cite{She95}. 
With the penalty parameter $\beta \gg 1$ we replace the incompressibility condition by 
\[
  \p = - \beta \ddiv \v. 
\]
In the semi-discrete case this corresponds to the constraint equation $M_p  p + B v = g$. It can be shown that this then leads to a DAE of (differentiation) index 1. However, this approach changes the solution of the system. In order to keep this difference of reasonable size, $\beta$ should be choosen relatively large. On the other hand, the condition number of the involved matrices increase with $\beta$ and thus, lead to numerical difficulties. The difficulty of a reasonable choice of $\beta$ is one drawback of this approach. 
A second disadvantage is that small velocities of order $\beta^{-1}$ or less cannot be resolved \cite{HeiV95}.
Further modifications of the penalty method, which are also applicable for slightly compressible fluids (Newtonian fluids)), are discussed, e.g., in \cite[Sect.~5]{HeiV95}. 
\subsubsection{Derivative of the constraint}\label{sect_solution_indexred_derivative}
Another very simple possibility to reduce the index of the given DAE is to replace the constraint by its derivative, the so-called {\em hidden constraint}. Instead of \eqref{eqn_NS_disc} we then consider the system 
\begin{subequations}
\label{eqn_derivativeConstraint}
\begin{align}
  M \dot {\tilde v} + K(\tilde v) - B^T \tilde p &= f, \\
  B \dot {\tilde v} &= \dot g. 
\end{align}
\end{subequations}
The initial condition remains unchanged, i.e., $\tilde v(0)=\qo$. It is well-known that this system has (differentiation) index 1. Nevertheless, numerical simulations, which rely on this formulation, show a linear drift from the solution manifold given by $Bv=g$. This can be seen as follows.

Although the two systems are equivalent, numerical errors are integrated over time and thus amplified. Solving the constraints only up to a small error, i.e.,
\[
  B v(t) = g(t) + \eps, \qquad
  B\dot{\tilde v}(t) = \dot g(t) + \eta   
\]
for small and constant $\eps$ and $\eta$, we calculate for $\tilde v$ that 
\[
  B\tilde v(t)   = B \qo + \int_0^t B\dot {\tilde v}(s) \ds
  = B \qo + \int_0^t \big(\dot g(s)  + \eta \big)\ds 
  = g(t) + t\,\eta. 
\]
Note that the last step holds because of the assumed consistency condition $B\qo = g(0)$. 
This shows that a constant error in the constraint of $\tilde v$ leads to an error which grows linearly in time. Because of this, the method of replacing the constraint by its derivative is -- although it is of index 1 -- not advisable. 
\subsubsection{Minimal extension}\label{sect_solution_indexred_minext}
Finally, we present the index reduction technique of {\em minimal extension} \cite[Ch.~6.4]{KunM06}. A general framework for an index reduction based on derivative arrays is given in~\cite[Ch.~6]{KunM06}; see also~\cite{Cam87}.  
Because of the special saddle point structure of system \eqref{eqn_NS_disc}, in which the constraint is explicitly given, this procedure can be simplified by using so-called {\em dummy variables}, cf.~\cite{MatS93, KunM04}.

Since we assume that $B$ is of full rank, there exists an invertible matrix $Q\in \R^{n,n}$ such that $BQ$ has the block structure $BQ = [ B_1,\  B_2 ]$ with an invertible matrix $B_2\in \R^{m,m}$. Note that the choice of $Q$ is not unique.  
We then use this transformation to partition the variable $v$, namely 
\[
  \begin{bmatrix} v_1 \\ v_2 \end{bmatrix} := Q^{-1} v,
\]
with according dimensions $v_1\in\R^{n-m}$ and $v_2\in\R^{m}$.  
With this, the constraint and its derivative may be written as 
\[
  B_2 v_2 = g - B_1v_1, \qquad
  B_2 \dot v_2 = \dot g - B_1 \dot v_1.
\]
Together with the differential equation \eqref{eqn_NS_disc_a} this yields an overdetermined system. 
Instead of an (expensive) search for projectors \cite[Ch.~6.2]{KunM06} or selectors \cite{Ste06}, which would bring the system back to its original size, we introduce a dummy variable $w_2 := \dot v_2$ to get rid of the redundancy. 
Note, however, that this leads to a slightly bigger system. 
More precisely, we extend the system dimensions from $n+m$ to $n+2m$, which is still moderate, since we usually have $m\ll n$. 
The extended system is again square and has the form 
\begin{align*}
  MQ \begin{bmatrix} \dot v_1\\ w_2\end{bmatrix} + K(v_1,v_2) - B^T p &= f,\\
  B_2 v_2 &= g - B_1v_1, \\
  B_2 w_2 &= \dot g - B_1 \dot v_1. 
\end{align*}
This DAE is of index $1$ and has the same solution set as the original system, cf.~\cite{AltH15}.

The drawback of this method is the need of a transformation matrix $Q$. 
With a suitable reordering of the variables, however, we can choose $Q$ to be the identity matrix. In this case, the needed variable transformation is just a permutation and thus, all variables keep their physical meaning. 
For some specific finite element schemes an algorithm to find such a permutation (which leads to an invertible $B_2$ block) is given in \cite{AltH15}. 
This paper also considers the implicit-explicit Euler scheme applied to the minimally extended system, which leads to a pair $(\calE, \mathcal A)$ of index 1. 
 \section{Numerical experiments}\label{sect_numerics}
To illustrate the performance and particular issues of the time-stepping schemes for the NSE, we consider the numerical simulation of the flow passing a cylinder in two space dimensions. This problem, also known as \emph{cylinder wake}, is a popular flow benchmark problem and a test field for flow control \cite{Wil96, NoaAMTT03, BenH15}. 

\begin{figure}
	\newlength\figureheight
	\newlength\figurewidth
	\setlength\figureheight{4cm}
	\setlength\figurewidth{15cm}
	\begin{tikzpicture}

\begin{axis}[
xmin=0, xmax=2.2,
ymin=0, ymax=.41,
width=\figurewidth,
height=\figureheight,
axis on top,
xtick={0,0.2, 2.2},
ytick={0,0.15, 0.25, .41},
tick label style={font=\footnotesize},
xlabel={$x_0$},
ylabel={$x_1$},
xlabel near ticks,
ylabel near ticks
]
\addplot graphics [includegraphics cmd=\pgfimage,xmin=0, xmax=2.2, ymin=0, ymax=0.41] {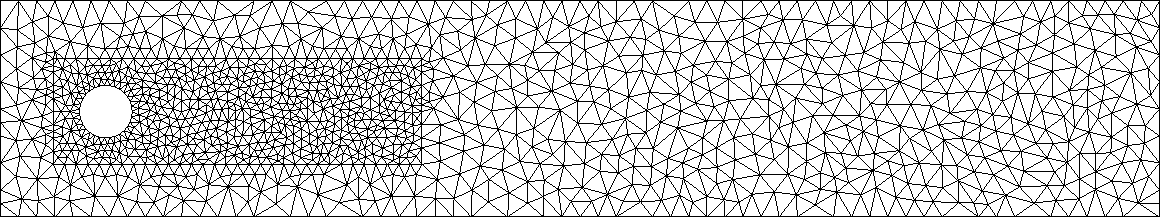};

\path [draw=black, fill opacity=0] (axis cs:13,1)--(axis cs:13,1);

\path [draw=black, fill opacity=0] (axis cs:0.05,13)--(axis cs:0.05,13);

\path [draw=black, fill opacity=0] (axis cs:13,1.38777878078145e-17)--(axis cs:13,1.38777878078145e-17);

\path [draw=black, fill opacity=0] (axis cs:0,13)--(axis cs:0,13);

\node at (axis cs:0.,0.25) [color=black, fill=white, anchor=south west] {$\Gammai$};
\node at (axis cs:2.2,0.21) [color=black, fill=white, anchor=south east] {$\Gammao$};

\end{axis}

\end{tikzpicture}
 	\caption{Illustration of the geometrical setup including the domains of distributed control and observation and of the velocity magnitude for the cylinder wake. Figure taken from \cite{BehBH17}.}
	\label{fig:cylwake}
\end{figure}
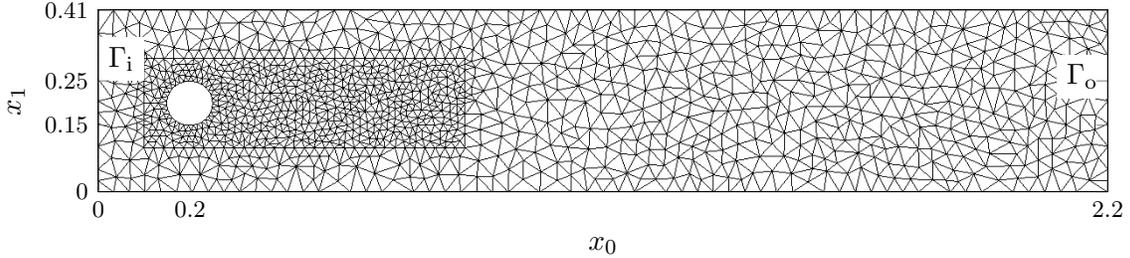

We consider the incompressible NSE \eqref{eqn_NSEnondim} on the domain as illustrated in Figure~\ref{fig:cylwake} with boundary $\Gamma$ and boundary conditions as follows. At the inflow $\Gammai$, we prescribe a parabolic velocity profile through the function 
$$g(s) = 4\bigl( 1 - \frac{s}{0.41}\bigr)\frac{s}{0.41} .$$
At the outflow $\Gammao$, we impose \emph{do-nothing} conditions, cf.~\eqref{eqn_bcs-donothing}. At the upper and the lower wall of the channel and at the cylinder periphery, we employ \emph{no-slip}, i.e., zero Dirichlet conditions. Thus, we set
\begin{equation*}
	\gamma (v, p) \colon
	\begin{cases}
		v = [g(x_1), 0 ]^T \quad &\text{on }\Gammai, \\
	p \vec n - \frac 1\RE \frac{\partial v}{\partial \vec n} =[0, 0]^T \quad &\text{on }\Gammao, \\
		v = [0, 0 ]^T \quad &\text{elsewhere on the boundary}.
	\end{cases}
\end{equation*}

We set $\RE=60$ and consider a $\mathcal P_2$--$\mathcal P_1$ (Taylor--Hood) finite element discretization of \eqref{eqn_NSEnondim} on the grid depicted in Figure~\ref{fig:cylwake} with $9356$ degrees of freedom in the velocity and $1289$ degrees of freedom in the pressure approximation. The result of the semi-discretization is a DAE of the form \eqref{eqn_NS_disc}, which we write as 
\begin{subequations}\label{eqn_numrs_NS_disc}
\begin{align}
  M \dot v + Av+N(v) - B^T p &= f, \\
  B v &= g
\end{align}
\end{subequations}
on the time interval $(0,1]$. 
Here, $A$ is the Laplacian or the linear part of $K$ as defined in \eqref{eqn_defdiscreteK} and $N$ denotes the convection part. The right-hand sides $f$ and $g$ account for the (static) boundary conditions. As initial value we take the corresponding steady-state \emph{Stokes} solution $\vS$, which is part of the solution to 
\begin{equation}
	\begin{bmatrix}
		A & -B^T \\ B & 0
	\end{bmatrix}
	\begin{bmatrix}
		\vS \\ \pS
	\end{bmatrix}
	=
	\begin{bmatrix}
		f \\ g
	\end{bmatrix}.
\end{equation}
For the schemes that need an initial value for the pressure, we provide it as $\pNS$ solving \eqref{eqn_numrs_NS_disc} at $t=0$ with $v(0)=\vS$ and $\dot v(0)=0$. Note that this gives a consistent initial pressure, since $g$ is constant, and thus, $B\dot v=0$.

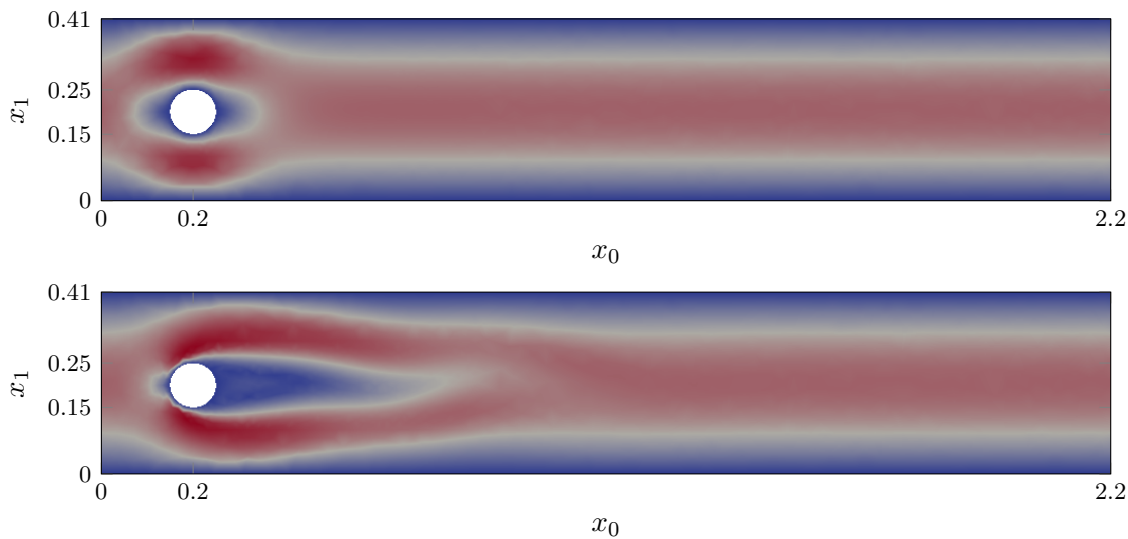
\begin{figure}
	\setlength\figureheight{4cm}
	\setlength\figurewidth{15cm}
	\begin{tikzpicture}

\begin{axis}[
xmin=0, xmax=2.2,
ymin=0, ymax=.41,
width=\figurewidth,
height=\figureheight,
axis on top,
xtick={0,0.2, 2.2},
ytick={0,0.15, 0.25, .41},
tick label style={font=\footnotesize},
xlabel={$x_0$},
ylabel={$x_1$},
xlabel near ticks,
ylabel near ticks
]
\addplot graphics [includegraphics cmd=\pgfimage,xmin=0, xmax=2.2, ymin=0, ymax=0.41] {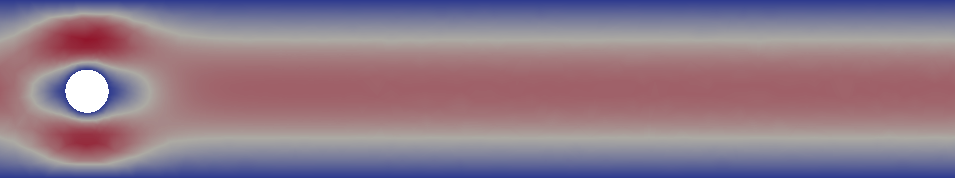};

\path [draw=black, fill opacity=0] (axis cs:13,1)--(axis cs:13,1);

\path [draw=black, fill opacity=0] (axis cs:0.05,13)--(axis cs:0.05,13);

\path [draw=black, fill opacity=0] (axis cs:13,1.38777878078145e-17)--(axis cs:13,1.38777878078145e-17);

\path [draw=black, fill opacity=0] (axis cs:0,13)--(axis cs:0,13);

\end{axis}

\end{tikzpicture}
 	\begin{tikzpicture}

\begin{axis}[
xmin=0, xmax=2.2,
ymin=0, ymax=.41,
width=\figurewidth,
height=\figureheight,
axis on top,
xtick={0,0.2, 2.2},
ytick={0,0.15, 0.25, .41},
tick label style={font=\footnotesize},
xlabel={$x_0$},
ylabel={$x_1$},
xlabel near ticks,
ylabel near ticks
]
\addplot graphics [includegraphics cmd=\pgfimage,xmin=0, xmax=2.2, ymin=0, ymax=0.41] {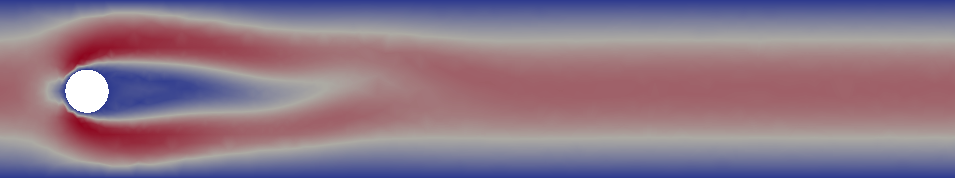};

\path [draw=black, fill opacity=0] (axis cs:13,1)--(axis cs:13,1);

\path [draw=black, fill opacity=0] (axis cs:0.05,13)--(axis cs:0.05,13);

\path [draw=black, fill opacity=0] (axis cs:13,1.38777878078145e-17)--(axis cs:13,1.38777878078145e-17);

\path [draw=black, fill opacity=0] (axis cs:0,13)--(axis cs:0,13);

\end{axis}

\end{tikzpicture}
 	\caption{Snapshots of the velocity magnitude computed with SIMPLE with $\tau = 1/1024$ and exact solves taken at $t=0$ (top) and $t=1$ (bottom).}
	\label{fig:cylwaket0t1}
\end{figure}

We consider the time-discretization of \eqref{eqn_numrs_NS_disc} by means of the implicit-explicit Euler scheme and compare it to the time-discretization via the SIMPLE scheme as described in Section~\ref{sec_SIMPLE}. Both schemes treat the nonlinearity explicitly. Comparing the computed approximations to a reference, obtained by the time-discretization via the implicit trapezoidal rule on a fine grid, we show that both schemes are convergent of order 1 as long as the resulting linear systems are solved with sufficient precision. For inexact solves, we show that the pressure approximation in the implicit-explicit Euler scheme diverges unlike for the SIMPLE approximation. For snapshots of the approximate velocity solution, see Figure~\ref{fig:cylwaket0t1}.

 The finite element implementation uses \emph{FEniCS, Version 2017.2} \cite{LogOeRW12}. For the iterative solutions of the linear system, we employ \href{https://github.com/andrenarchy/krypy}{\emph{Krypy}} \cite{Krypy17}.
The code used for the numerical investigations is freely available for reproducing the reported results and as a benchmark for further developments in the time integration of semi-explicit DAEs of strangeness-index 1; see Figure \ref{fig:linkcodndat} for a stable link to the online repository.

\begin{figure}[h!]
  \begin{framed}
    \textbf{Code and Data Availability} \\
    The source code and the data of the implementations used to compute the presented results is available from:
    \begin{center}
			\href{https://doi.org/10.5281/zenodo.998909}{\texttt{doi:10.5281/zenodo.998909}}
    \end{center}
	contact the author Jan Heiland for licensing information.
  \end{framed}
	\caption{Link to code and data.}\label{fig:linkcodndat}
\end{figure}
\subsection{Time integration with implicit-explicit Euler}\label{sect_numerics_imexeul}
With the implicit-explicit Euler time discretization, at every time step the linear system 
\begin{equation}\label{eqn_imexEulUpd}
	 \begin{bmatrix} \tfrac 1\tau M+A & B^T \\ B & 0 \end{bmatrix}
		 \begin{bmatrix} v^{k+1} \\ -p^k \end{bmatrix} = \rhs:=
	       \begin{bmatrix} \tfrac 1\tau M - N(v^k ) +f \\ g \end{bmatrix}
\end{equation}
is solved for the velocity and the pressure approximations. For the approximate solution of the linear systems, we use \emph{MinRes} iterations, which are stopped as soon as the relative residual drops below a given tolerance \lintol, i.e., 
$$
	 \norm*{\begin{bmatrix} \tfrac 1\tau M+A & B^T \\ B & 0 \end{bmatrix}
		 \begin{bmatrix} v^{k+1} \\ -p^k \end{bmatrix} - \rhs\,}_{(M^{-1}, M_{\mathcal Q}^{-1})} \leq \lintol~ \norm{\rhs}_{(M^{-1}, M_{\mathcal Q}^{-1})},
$$
where we use the norm induced by the inverses of the mass matrix $M$ and of the mass matrix of the pressure approximation space $M_{\mathcal Q}$.

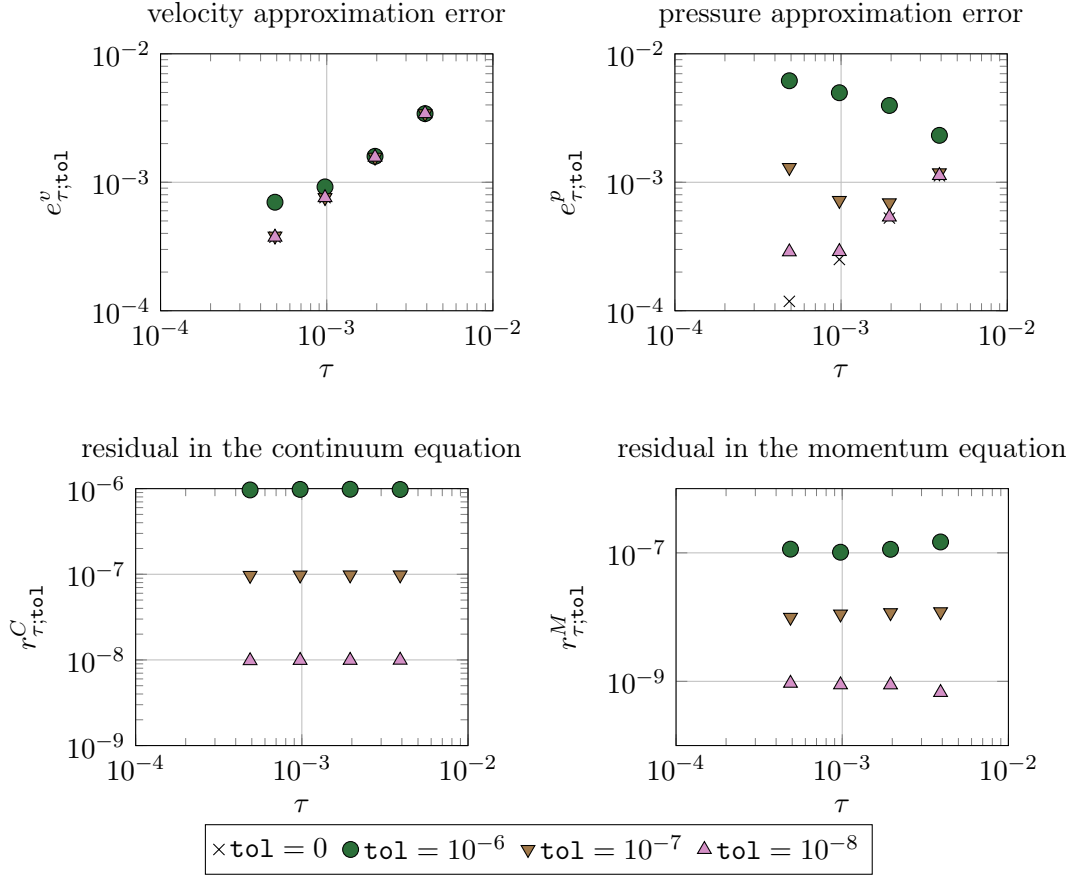
\begin{figure}
	\setlength\figureheight{5cm}
	\setlength\figurewidth{6cm}
	\begin{tikzpicture}

\begin{axis}[
	legend columns = 4,
	legend entries={$\lintol= 0~~~$,$\lintol= 10^{-6}~~~$,$\lintol= 10^{-7}~~~$,$\lintol= 10^{-8}$},
	legend to name=imexeulleg,
title={velocity approximation error},
xmin=0.0001, xmax=0.01,
ymin=0.0001, ymax=0.01,
xmode=log,
ymode=log,
axis on top,
xlabel=$\tau$,
ylabel=$\ver$,
width=\figurewidth,
height=\figureheight,
xmajorgrids,
ymajorgrids
]
\addplot [black, mark=x, mark size=3, only marks]
table {0.00390625 0.00341974712593372
0.001953125 0.00156723623841947
0.0009765625 0.00075930797046964
0.00048828125 0.000374001701706994
};
\addplot [color1, mark=*, mark size=3, mark options={draw=black}, only marks]
table {0.00390625 0.00342259300619518
0.001953125 0.00159263112673649
0.0009765625 0.000920545773814314
0.00048828125 0.00069941031152655
};
\addplot [color2, mark=triangle*, mark size=3, mark options={rotate=180,draw=black}, only marks]
table {0.00390625 0.00341944313711747
0.001953125 0.00156734641517254
0.0009765625 0.000759992829778325
0.00048828125 0.000381027777509547
};
\addplot [color3, mark=triangle*, mark size=3, mark options={draw=black}, only marks]
table {0.00390625 0.0034197256455773
0.001953125 0.00156715505035038
0.0009765625 0.000759307614053305
0.00048828125 0.000374008790490444
};
\end{axis}

\end{tikzpicture}
 	\begin{tikzpicture}

\definecolor{color0}{rgb}{0,0.75,0.75}

\begin{axis}[
title={pressure approximation error},
xmin=0.0001, xmax=0.01,
ymin=0.0001, ymax=0.01,
xmode=log,
ymode=log,
axis on top,
xlabel=$\tau$,
ylabel=$\per$,
width=\figurewidth,
height=\figureheight,
xmajorgrids,
ymajorgrids
]
\addplot [black, mark=x, mark size=3, only marks]
table {0.00390625 0.00112342607803267
0.001953125 0.000526725389693455
0.0009765625 0.000250583557066635
0.00048828125 0.000118332104058548
};
\addplot [color1, mark=*, mark size=3, mark options={draw=black}, only marks]
table {0.00390625 0.00231902326587635
0.001953125 0.00396215431974663
0.0009765625 0.00497704319281976
0.00048828125 0.00616932481969313
};
\addplot [color2, mark=triangle*, mark size=3, mark options={rotate=180,draw=black}, only marks]
table {0.00390625 0.00118494598112798
0.001953125 0.0006920720044629
0.0009765625 0.000722122295107473
0.00048828125 0.00130282298313308
};
\addplot [color3, mark=triangle*, mark size=3, mark options={draw=black}, only marks]
table {0.00390625 0.00112817647233799
0.001953125 0.000536549225451843
0.0009765625 0.000288921186799507
0.00048828125 0.000287446958053257
};
\end{axis}

\end{tikzpicture}
  \\[.2in]
	\begin{tikzpicture}

\begin{axis}[
title={residual in the continuum equation},
xmin=0.0001, xmax=0.01,
ymin=1e-09, ymax=1e-06,
xlabel=$\tau$,
ylabel=$\rescon$,
xmode=log,
ymode=log,
axis on top,
width=\figurewidth,
height=\figureheight,
xmajorgrids,
ymajorgrids
]
\addplot [color1, mark=*, mark size=3, mark options={draw=black}, only marks]
table {0.00390625 9.77192080248568e-07
0.001953125 9.82807817109837e-07
0.0009765625 9.79453077561229e-07
0.00048828125 9.65343982257065e-07
};
\addplot [color2, mark=triangle*, mark size=3, mark options={rotate=180,draw=black}, only marks]
table {0.00390625 9.82753384588453e-08
0.001953125 9.82290526617209e-08
0.0009765625 9.79299937886903e-08
0.00048828125 9.70468920089708e-08
};
\addplot [color3, mark=triangle*, mark size=3, mark options={draw=black}, only marks]
table {0.00390625 9.92699240838445e-09
0.001953125 9.88571561477018e-09
0.0009765625 9.86746207252476e-09
0.00048828125 9.79751014316184e-09
};
\end{axis}

\end{tikzpicture}
 	\begin{tikzpicture}

\begin{axis}[
title={residual in the momentum equation},
xmin=0.0001, xmax=0.01,
ymin=1e-10, ymax=1e-06,
xmode=log,
ymode=log,
axis on top,
width=\figurewidth,
height=\figureheight,
xlabel=$\tau$,
ylabel=$\resmom$,
xmajorgrids,
ymajorgrids
]
\addplot [color1, mark=*, mark size=3, mark options={draw=black}, only marks]
table {0.00390625 1.4783262317102e-07
0.001953125 1.13554134841253e-07
0.0009765625 1.02317365544488e-07
0.00048828125 1.14192008650874e-07
};
\addplot [color2, mark=triangle*, mark size=3, mark options={rotate=180,draw=black}, only marks]
table {0.00390625 1.21128841945772e-08
0.001953125 1.17121381811043e-08
0.0009765625 1.1155132775438e-08
0.00048828125 9.91638912927493e-09
};
\addplot [color3, mark=triangle*, mark size=3, mark options={draw=black}, only marks]
table {0.00390625 6.70607971168895e-10
0.001953125 8.804522602121e-10
0.0009765625 8.82558339939375e-10
0.00048828125 9.33185099102533e-10
};
\end{axis}

\end{tikzpicture}
 	\ref{imexeulleg}
	\caption{Error in the velocity and the pressure approximation provided by the implicit-explicit Euler algorithm for iterative solves with varying tolerances. The crosses are the errors obtained with direct solves.}
	\label{fig:imexeuler_errs}
\end{figure}

It can be seen from the error plots in Figure~\ref{fig:imexeuler_errs}, for exact solves, that the implicit-explicit Euler converges linearly in the velocity and the pressure approximation, which is in line with the theory~\cite[Tab.~2.3]{HaiLR89}. For inexact solves, however, for smaller time-steps, the pressure approximation diverges linearly. As we have shown in \cite{AltH15} this is an inherent instability of the index-2 formulation.

The reported residuals and errors are defined as follows:
\begin{equation*}
	\ver = \trapz(\norm{v_{\tau;\lintol}-v_\rfrc}_M)
\end{equation*}
where the subscript \rfrc~denotes the reference solution, the subscripts $\tau;\lintol$ denote the approximation that is computed on the grid of size $\tau$ with the linear equations solved with tolerance $\lintol$, and where $\trapz(s)$ denotes the approximation to the integral $\int_0^{t_N} s(t)\inva t$ by means of the piecewise trapezoidal rule with step size $\tau$. Analogously, we define the error in the pressure approximation
$$
\per = \trapz(\norm{p_{\tau;\lintol}-p_\rfrc}_{M_P})
$$
and the integrated residuals in the momentum equation $\resmom$ as the integral of the function
$$
t^k \mapsto \norm{(\tfrac  1\tau M+A)v_{\tau;\lintol}^{k+1} - \tfrac 1\tau Mv_{\tau;\lintol}^k - B^T p_{\tau;\lintol}^k + N(v^k_{\tau;\lintol}) - f^k}_{M^{-1}}
$$
and in the continuum equation as
$$
\rescon := \trapz(\norm{Bv_{\tau;\lintol} -g }_{M_{\mathcal Q}^{-1}}).
$$
\subsection{Time integration with SIMPLE}\label{sect_numerics_SIMPLE}
For this particular time discretization, we have to solve three linear systems, namely
\begin{subequations}
\begin{align}
   	(\tfrac  1\tau M+A)\tilde v^{k+1} &= \rhs_1:=\tfrac 1\tau Mv^k + B^T p^k + f^k, \\
		B(\tfrac 1\tau M+A)^{-1}B^Tp_\Delta^{k+1} &= \rhs_2 :=-B\tilde v^{k+1} + g^{k+1}, \label{eqn_upd_SIMPLE_pd}\\
		\intertext{and}
		(\tfrac  1\tau M+A)v_\Delta^{k+1}  &= \rhs_3:= B^Tp_\Delta^{k+1} 
\end{align}
\end{subequations}
to compute the updates as $v^{k+1} = \tilde v^{k+1} + v_\Delta^{k+1} $ and $p^{k+1} = p^k  + p_\Delta^{k+1} $. For the approximate solution, we use \emph{CG} iterations until the relative residuals, measured in the $M^{-1}$ norm (or $M_{\mathcal Q}^{-1}$ for \eqref{eqn_upd_SIMPLE_pd}), drop below a given tolerance $\lintol$. 

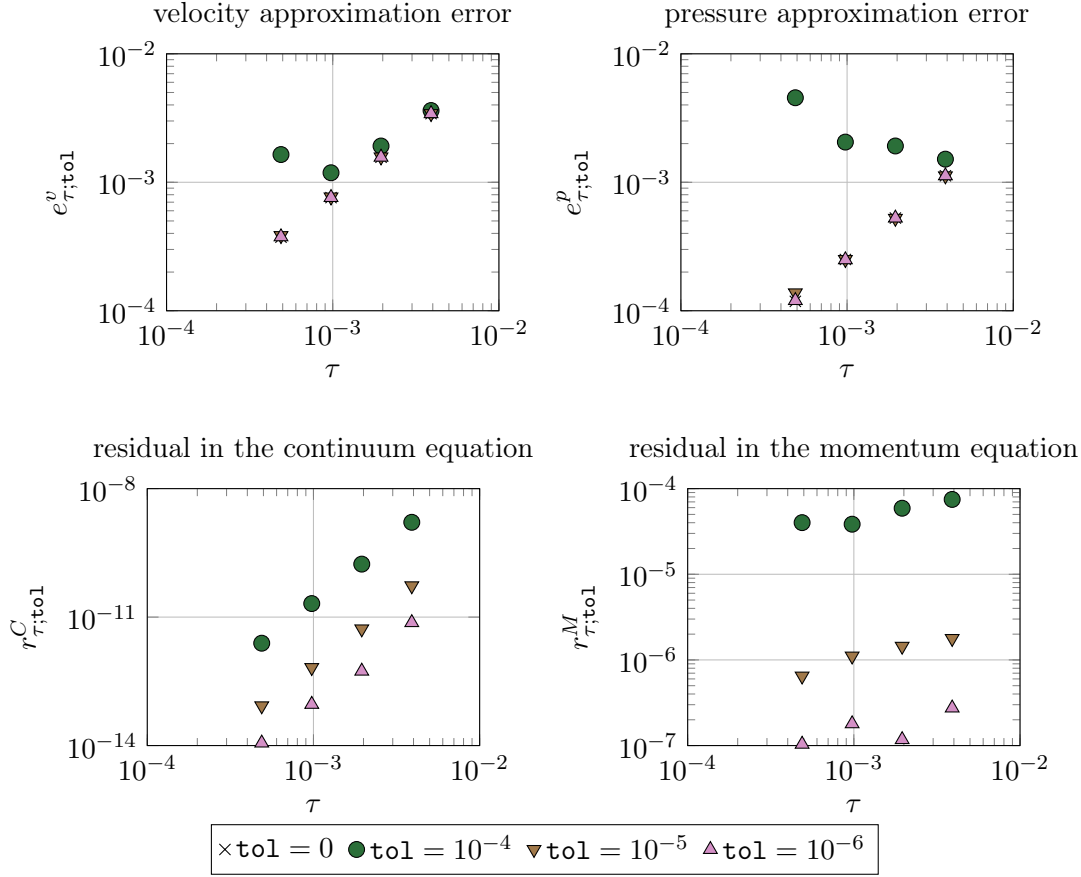
\begin{figure}
	\setlength\figureheight{5cm}
	\setlength\figurewidth{6cm}
	\begin{tikzpicture}

\begin{axis}[
legend columns = 4,
legend entries={$\lintol= 0~~~$,$\lintol= 10^{-4}~~~$,$\lintol= 10^{-5}~~~$,$\lintol= 10^{-6}$},
legend to name=simpleleg,
title={velocity approximation error},
xmin=0.0001, xmax=0.01,
ymin=0.0001, ymax=0.01,
xlabel=$\tau$,
ylabel=$\ver$,
xmode=log,
ymode=log,
axis on top,
width=\figurewidth,
height=\figureheight,
xmajorgrids,
ymajorgrids
]
\addplot [black, mark=x, mark size=3, only marks]
table {0.00390625 0.00341974712593376
0.001953125 0.00156723623841913
0.0009765625 0.000759307970469215
0.00048828125 0.000374001701707145
};
\addplot [color1, mark=*, mark size=3, mark options={draw=black}, only marks]
table {0.00390625 0.003618273187754
0.001953125 0.00191569404087298
0.0009765625 0.00118583301636522
0.00048828125 0.00164384295116239
};
\addplot [color2, mark=triangle*, mark size=3, mark options={rotate=180,draw=black}, only marks]
table {0.00390625 0.00342497727472714
0.001953125 0.00157379543589841
0.0009765625 0.000768204203986231
0.00048828125 0.000384533529476464
};
\addplot [color3, mark=triangle*, mark size=3, mark options={draw=black}, only marks]
table {0.00390625 0.00342059555925665
0.001953125 0.00156800205010477
0.0009765625 0.000760852130246938
0.00048828125 0.000376358385119053
};
\end{axis}

\end{tikzpicture}
 	\begin{tikzpicture}

\begin{axis}[
title={pressure approximation error},
xmin=0.0001, xmax=0.01,
ymin=0.0001, ymax=0.01,
xmode=log,
ymode=log,
axis on top,
width=\figurewidth,
height=\figureheight,
xlabel=$\tau$,
ylabel=$\per$,
xmajorgrids,
ymajorgrids
]
\addplot [black, mark=x, mark size=3, only marks]
table {0.00390625 0.00112342607803244
0.001953125 0.000526725389692994
0.0009765625 0.00025058355706642
0.00048828125 0.000118332104059264
};
\addplot [color1, mark=*, mark size=3, mark options={draw=black}, only marks]
table {0.00390625 0.00151352634584607
0.001953125 0.00191692954146792
0.0009765625 0.00205476586460881
0.00048828125 0.00454949213874456
};
\addplot [color2, mark=triangle*, mark size=3, mark options={rotate=180,draw=black}, only marks]
table {0.00390625 0.00112137804119208
0.001953125 0.000523750680621591
0.0009765625 0.000250418706649227
0.00048828125 0.000137730983012303
};
\addplot [color3, mark=triangle*, mark size=3, mark options={draw=black}, only marks]
table {0.00390625 0.00112311692046461
0.001953125 0.000526590946001437
0.0009765625 0.00024934029256772
0.00048828125 0.000120039222856922
};
\end{axis}

\end{tikzpicture}
 \\[.2in]
	\begin{tikzpicture}

\begin{axis}[
title={residual in the continuum equation},
xmin=0.0001, xmax=0.01,
ymin=1e-14, ymax=1e-08,
xmode=log,
ymode=log,
xlabel=$\tau$,
ylabel=$\rescon$,
axis on top,
width=\figurewidth,
height=\figureheight,
xmajorgrids,
ymajorgrids
]
\addplot [color1, mark=*, mark size=3, mark options={draw=black}, only marks]
table {0.00390625 1.63587321372335e-09
0.001953125 1.72598677992039e-10
0.0009765625 2.0693198744439e-11
0.00048828125 2.4559609967083e-12
};
\addplot [color2, mark=triangle*, mark size=3, mark options={rotate=180,draw=black}, only marks]
table {0.00390625 5.37740608304353e-11
0.001953125 5.29423905548181e-12
0.0009765625 6.62764492128139e-13
0.00048828125 8.46452579675574e-14
};
\addplot [color3, mark=triangle*, mark size=3, mark options={draw=black}, only marks]
table {0.00390625 7.43333164057224e-12
0.001953125 5.47148101801832e-13
0.0009765625 9.17851212440721e-14
0.00048828125 1.14594465006747e-14
};
\end{axis}

\end{tikzpicture}
 	\begin{tikzpicture}

\begin{axis}[
title={residual in the momentum equation},
xmin=0.0001, xmax=0.01,
ymin=1e-07, ymax=0.0001,
xmode=log,
ymode=log,
axis on top,
width=\figurewidth,
height=\figureheight,
xlabel=$\tau$,
ylabel=$\resmom$,
xmajorgrids,
ymajorgrids
]
\addplot [color1, mark=*, mark size=3, mark options={draw=black}, only marks]
table {0.00390625 7.48015479335723e-05
0.001953125 5.90364736146086e-05
0.0009765625 3.84264982813972e-05
0.00048828125 4.0080482980536e-05
};
\addplot [color2, mark=triangle*, mark size=3, mark options={rotate=180,draw=black}, only marks]
table {0.00390625 1.77941900186072e-06
0.001953125 1.44295683392865e-06
0.0009765625 1.10894688200532e-06
0.00048828125 6.47425490233302e-07
};
\addplot [color3, mark=triangle*, mark size=3, mark options={draw=black}, only marks]
table {0.00390625 2.7459911087053e-07
0.001953125 1.17353253120218e-07
0.0009765625 1.79701506965257e-07
0.00048828125 1.03666872168552e-07
};
\end{axis}

\end{tikzpicture}
 	\ref{simpleleg}
	\caption{Error in the velocity and the pressure approximation provided by the SIMPLE algorithm for exact solves and iterative solves with varying tolerances \lintol.}
	\label{fig:SIMPLE_errs}
\end{figure}

As one can see from the error plots in Figure~\ref{fig:SIMPLE_errs}, the inexact solves affect the approximation only for a rough tolerance $\lintol = 10^{-4}$, which is only one order of magnitude smaller than the actual approximation error. Thus the breakdown in the convergence observed in Figure \ref{fig:SIMPLE_errs} is probably due to the accumulation of errors that every single step scheme suffers from. For smaller tolerances, the approximations almost achieve the accuracy of the direct solves. Interestingly, the continuity equation \eqref{eqn_NS_disc_b}, which is only an implicit part of the SIMPLE scheme, is fulfilled at much better accuracy than the choices of the tolerances suggest; see Figure~\ref{fig:SIMPLE_errs}.

 \section{Conclusion}\label{sect_conclusion}
In this paper, we have discussed the incompressible NSE from a DAE point of view, in which the incompressibility is interpreted as an algebraic constraint. Thus, a spatial discretization leads to a DAE. If the time is discretized as well, a \dDAE~is obtained -- a sequence of equations, that define the numerical approximations. 

We have discussed suitable approaches to well-posed semi-discrete approximations and investigated the \dDAE~stemming from different time-discretizations in terms of the Kronecker index. It turned out that commonly and successfully used time integration schemes like the SIMPLE algorithm define a \dDAE~of index 1, whereas a time-discretization of the semi-discrete NSE by an implicit-explicit Euler scheme leads to a \dDAE~of index 2. Alternatively, one may first apply an index reduction on the semi-discrete level and then apply time marching schemes.

The advantage of the discrete index-1 formulations is that they avoid implicit derivations that amplify computational errors. This mechanism will likely lead to larger errors, in particular if the equations are solved with limited accuracy. We have illustrated the origin of this behavior in a small analytical example and verified it in a numerical simulation. The code of the numerical test case, an implementation of the 2D cylinder wake, is provided to serve as a benchmark for future developments of time integration schemes for DAEs of Navier--Stokes type.
\newcommand{\etalchar}[1]{$^{#1}$}

\appendix
\section{Strangeness index of equation \eqref{eqn_NS_disc}}\label{app_sindex}
We analyse in detail the strangeness-index of the DAE \eqref{eqn_NS_disc}. 
Note that we do not ask for any assumptions on the nonlinearity $K$ and that we allow the matrix $B$ to be rank-deficient. 
This then also implies that the d-index of \eqref{eqn_NS_disc} equals 2 if it is well-defined, i.e., if $B$ is of full rank.
\subsection{Linear case}
Considering any linearization of the Navier-Stokes equations, i.e., $K(u) = Ku$ in \eqref{eqn_NS_disc}, we deal with the matrix pair
\[
  (E,A) = \Big( \begin{bmatrix} M & 0 \\ 0 & 0 \end{bmatrix},  \begin{bmatrix} K & B^T \\ B & 0 \end{bmatrix} \Big).
\]
Following \cite[Th.~3.11]{KunM06}, we can construct a to $(E,A)$ (globally) equivalent pair $(\tilde E, \tilde A )$ of the form 
\[
  \tilde E =  \left[ \begin{array}{cc|cc} I_b & & \\ & I_{n-b} & \\ \hline & & 0 \\ & & 0 \end{array} \right], \qquad
  \tilde A =  \left[ \begin{array}{cc|cc} 0 & A_{12} & A_{13} \\ 0 & 0 & A_{23} \\ \hline I_b & & 0 \\ & 0 & 0 \end{array} \right]
\]
with $A_{13} \in \R^{b,m}$ being of full rank. 
Thus, the original system \eqref{eqn_NS_disc} is equivalent to a system of the form 
\[
  \dot x_1 = A_{12}x_2 + A_{13}x_3 + f_1, \qquad 
  \dot x_2 = A_{23}x_3 + f_2, \qquad
  0 = x_1 + f_3, \qquad
  0 = f_4
\]
with dimensions $x_1(t)\in\R^{b}$, $x_2(t)\in\R^{n-b}$, $x_3(t)\in\R^{m}$. 
Since we have a differential and an algebraic equation for $x_1$ (this causes the 'strangeness'), we use the derivative of $0 = x_1 + f_3$ in order to eliminate $\dot x_1$ in the first equation. 
Hence, we consider the pair $(E_\text{mod},A_\text{mod})$ with
\[
  E_\text{mod} =  \left[ \begin{array}{cc|cc} 0 & & \\ & I_{n-b} & \\ \hline & & 0 \\ & & 0 \end{array} \right], \qquad
  A_\text{mod} =  \left[ \begin{array}{cc|cc} 0 & A_{12} & A_{13} \\ 0 & 0 & A_{23} \\ \hline I_b & & 0 \\ & 0 & 0 \end{array} \right].
\]
Since $A_{13}$ is of full rank, one can show that system $(E_\text{mod},A_\text{mod})$ is strangeness-free, cf.~ the calculation in \cite[Th.~3.7]{KunM06}. 
Since we have obtained a strangeness-free system with only one differentiation, system \eqref{eqn_NS_disc} has strangeness-index one. 
\subsection{Nonlinear case}
The general form of a nonlinear DAE is given by 
\[
  F(t, x, \dot x) = 0.
\]
In regard of system \eqref{eqn_NS_disc} we set $x:=[q^T,\ p^T]^T$ and define 
\[
  F(t, x, \dot x) 
  := \begin{bmatrix} M \dot q - K(q) - B^Tp - f  \\ - Bq + g \end{bmatrix}
  = E\dot x - A(x) - h
\]
with 
\[
  E := \begin{bmatrix} M & 0 \\ 0 & 0 \end{bmatrix},  \qquad
  A_x := \frac{\partial A(x)}{\partial x} = \begin{bmatrix} K_q & B^T \\ B & 0 \end{bmatrix}. 
\]
In the sequel we show that \eqref{eqn_NS_disc} has strangeness-index 1 also in the nonlinear case. 
For this, we assume that $B$ has full rank such that there are no vanishing equations and the pressure variable is uniquely defined. 
In the case $\rank B = b < m$, we consider the following transformation. 

Let $C_0\in \R^{m, m-b}$ be the matrix of full rank satisfying $B^TC_0 = 0$. 
Furthermore, $C'\in \R^{m, b}$ defines any matrix such that $C = [C_0\ C'] \in \R^{m, m}$ is invertible. 
With this, we obtain the relation 
\[
  B^TC 
  = \begin{bmatrix} B^T C_0 & B^T C' \end{bmatrix}
  = \begin{bmatrix} 0 & \tilde B^T \end{bmatrix}
\]
with $\tilde B \in \R^{b, n}$ having full rank. 
With the matrix $C$ in hand, we first introduce the new pressure variable $\tilde p := C^{-1}p$. 
Thus, we consider the the pair $z := [q^T, \tilde p^T]^T$. 
As a second step, we multiply equation \eqref{eqn_NS_disc} by the block-diagonal matrix $\diag(I_n, C^T)$ from the left. 
In total, this yields the equivalent DAE
\[
  M \dot q = K(q) + \begin{bmatrix}  0 & \tilde B^T \end{bmatrix} \tilde p + f, \qquad
  \begin{bmatrix} 0 \\ \tilde B \end{bmatrix} q = C^Tg.
\]
Note that the constraint contains $(m-b)$ consistency equations of the form $0=g_1$. 
Assuming that system \eqref{eqn_NS_disc} is solvable, we suppose that these are in fact vanishing equations. 
Thus, they have no influence on the index of the system. 
Furthermore, the first $(m-b)$ components of the transformed pressure $\tilde p$ do not influence the system. These components are underdetermined and may be omitted, again without changing the index. 
Leaving out the underdetermined parts as well as the vanishing equations, we obtain a system of the form \eqref{eqn_NS_disc} with a full rank matrix $B$. 

In the sequel, we assume that $\rank B = m$ and show that \cite[Hyp.~4.2]{KunM06} is satisfied for $\mu=1$. 
Note that this hypothesis is not satisfied for $\mu=0$, i.e., the system is not strangeness-free. 
We define the matrices
\[
  M_1 := \begin{bmatrix} E & 0 \\ -A_x & E \end{bmatrix}, \qquad
  N_1 := \begin{bmatrix} A_x & 0 \\ 0 & 0 \end{bmatrix}. 
\]
We now pass through the list of points of the hypothesis in \cite[Hyp.~4.2]{KunM06}:
\begin{enumerate}
  \item First, we note that the rank of $M_1$ equals $2n$ and we set $a:=2(n+m)-2n=2m$. Thus, the system contains $2m$ algebraic variables (the pressure and the part of $q$, which is not divergence-free). Furthermore, we define $Z_2\in\R^{2(n+m),2m}$ by $Z_2^TM_1 = 0$, i.e.,
\[
  Z_2 = \left[ \begin{array}{c|c} 0 & M^{-1}B^T \\ I_m & 0 \\ \hline 0 & 0 \\ 0 & I_m \end{array} \right]. 
\]
  \item As a second step we define $\hat A_2:= Z_2^T N_1 [I_{n+m},\ 0]^T$, which yields
\[ 
  \hat A_2 
  = \left[ \begin{array}{cc|cc} 0 & I_m & 0 & 0 \\ \hline B M^{-1} & 0 & I_m & 0 \end{array} \right] 
    \left[ \begin{array}{cc} K_q & B^T \\ B & 0 \\ \hline 0 & 0 \\ 0 & 0  \end{array} \right] 
  = \left[ \begin{array}{c|c} B & 0 \\ \hline B M^{-1}K_u & B M^{-1}B^T \end{array} \right]. 
\]
  This matrix has rank $2m$, since the full-rank property of $B$ implies that $BM^{-1}B^T$ is invertible. 
  We define $d:=n-m$ as the number of differential variables and $T_2\in\R^{n+m,n-m}$ by $\hat A_2T_2 = 0$. 
  Let $C\in \R^{n,n-m}$ be a matrix of full rank with $BC=0$ and $C_2 := -(BM^{-1}B^T)^{-1}BM^{-1}K_uC \in \R^{m,n-m}$. 
  Then, we set
\[
  T_2 := \begin{bmatrix} C \\ C_2 \end{bmatrix}.  
\]
  \item Finally, we compute the rank of $ET_2$. Since $C$ has full rank, this equals $\rank MC = n-m = d$. The matrix $Z_1^T := [C^T\ 0] \in \R^{n-m, n+m}$ satisfies 
\[
  \rank Z_1^T E T_2 = \rank C^TMC = n-m = d.   
\] 
\end{enumerate}
Thus, the hypothesis in \cite[Hyp.~4.2]{KunM06} is satisfied for $\mu=1$, which implies that the nonlinear DAE \eqref{eqn_NS_disc} has strangeness-index one.  \section{Difference-algebraic equation index of the considered systems}\label{app_diffAE}
In this appendix, we derive the Kronecker index for the discrete schemes considered in Section~
\ref{sect_solution_methods}.
\subsection{IMEX Euler}
We start with the implicit-explicit Euler discretization, that gives a scheme $\mathcal Ex^{k+1} = \mathcal A^k x^k+ h^k $ with the matrix pair 
\begin{equation*}
	\EApair = 
	\bigl( \begin{bmatrix} \tfrac 1\tau M & 0 \\ 0 & 0 \end{bmatrix},  
	       \begin{bmatrix} \tfrac 1\tau M+A & B^T \\ B & 0 \end{bmatrix} \bigr)
\end{equation*}
as in \eqref{eqn_EAplainEul}. For sufficiently small $\tau$, due to the definiteness of $M$ and the full-rank property of $B$, the matrix $\mathcal A$ is invertible and thus, the pair $(\mathcal E, \mathcal A)$ is regular. 
Let $S$ denote the matrix $BM^{-1}B^T$. 
If one applies 
\begin{equation*}
	\begin{bmatrix} M^{-\frac 12} & B^T S \\ 0 & I \end{bmatrix} 
		\begin{bmatrix} I & 0 \\ BM^{-1} & I \end{bmatrix} 
			\rightarrow \EApair \leftarrow
			\begin{bmatrix} M^{-\frac 12} & 0 \\ -S^{-1}BM^{-1}(\tfrac 1\tau M + A) & I \end{bmatrix} 
\end{equation*}
from the left and the right, one finds that $(\mathcal E, \mathcal A)$ is similar to
\begin{align*}
	&\bigl (
	\begin{bmatrix} \frac 1\tau (I-M^{-\frac 12}B^T S BM^{-\frac 12}) & 0 \\ \tfrac 1\tau B & 0\end{bmatrix} , \\
     &\qquad\qquad\begin{bmatrix} (I-M^{-\frac 12}B^T S BM^{-\frac 12})(\tfrac 1\tau I + M^{-\frac 12}AM^{-\frac 12})-M^{-\frac 12}B^T S BM^{-\frac 12} & 0 \\ 0 & S \end{bmatrix} \bigr).
\end{align*}
Since $B$ is of full rank, there exists an orthogonal matrix $Q$ and an invertible matrix $R$ such that $BM^{-\frac 12}Q = \big[ 0\ \   R \big]$ and, in particular, 
$$
		Q^T(I-M^{-\frac 12}B^T S BM^{-\frac 12})Q = 
		\begin{bmatrix} I & 0 \\ 0 & 0 \end{bmatrix}. 
$$
Thus, the corresponding similarity transformation transforms  $(\mathcal E, \mathcal A)$ into
		\begin{equation*}
			\bigl (
			\begin{bmatrix}
				\tfrac 1\tau I & 0 & 0\\
				0 & 0& 0\\
				0 & \tfrac 1\tau R & 0 
			\end{bmatrix},
			\begin{bmatrix}
				* & 0 & 0\\
				*_2 & I & 0\\
				0 & 0 & S 
			\end{bmatrix}
			\bigr),
		\end{equation*}
		where $*$ stands for the block matrix entries that need not be specified further. With another few regular row and column transformations, one can eliminate the entry $*_2$ and read off the \emph{Kronecker} index of \EApair~as the index of nilpotency of $\begin{bmatrix} 0& 0\\ \tfrac 1\tau R & 0 \end{bmatrix}$ which is $2$.
\subsection{Projection scheme}\label{sec_appdx_kronind_prjctn}
The matrix coefficient pair of the Projection scheme  
\eqref{eqn_eapair_prjctn}
reads
\begin{equation}
	\bigl(\mathcal E, \mathcal A\bigr) = 
	\bigl( 
	\begin{bmatrix} \tfrac 1\tau M & 0 & 0 & 0\\
									0 & 0 & 0 & 0 \\
									-M & -\frac \tau 2 B^T & M & 0 \\
									0 & -I & 0 & I

	\end{bmatrix},  
	\begin{bmatrix} 0 & \tfrac 1\tau M+A & 0 & B^T \\
									-\frac 2\tau B & -BM^{-1}B^T & 0 & 0 \\
									0 & 0 & 0 & 0 \\
									0 & 0 & 0 & I
	\end{bmatrix}
	\bigr).
\end{equation}
If we define $S := BM^{-1}B^T$, if we move the second row and column to the left and bottom, respectively, and if we rescale certain rows and columns, we find that the pair is equivalent to
\begin{equation*}
	\bigl(\mathcal E, \mathcal A\bigr) \sim 
	\bigl( 
	\begin{bmatrix} I & 0 & 0 & 0\\
									* & I & 0 & *_1  \\
									0 & 0 & I & *_2 \\
									0 & 0 & 0 & 0 

	\end{bmatrix},  
	\begin{bmatrix} 0 & 0 & * & *_3 \\
									0 & 0 & 0 & 0  \\
									0 & 0 & * & 0  \\
								  *_4 & 0 & 0 & -S 
	\end{bmatrix}
	\bigr),
\end{equation*}
where the $*$'s stand for unspecified but possibly nonzero entries. Since, in particular, $S$ is invertible, one can eliminate the entries $*_1$--$*_4$ by regular row and column manipulations without affecting the invertibility of the left upper $3\times3$ block in the transformed $\mathcal E$ and read off the Kronecker index of \eqref{eqn_eapair_prjctn} as the index of nilpotency of $0$ which is $1$.
\subsection{SIMPLE}\label{sec_appdx_kronind_SMPL}
The matrix coefficient pair of the SIMPLE scheme \eqref{eqn_eapair_SMPL} reads 
\begin{equation*}
	\bigl(\mathcal E, \mathcal A\bigr) = 
	\bigl( 
	\begin{bmatrix} \tfrac 1\tau M & 0 & 0 & 0\\
									0 & 0 & 0 & 0 \\
									I & (\tfrac 1\tau M+A)^{-1}B^T & -I & 0 \\
									0 & -I & 0 & I

	\end{bmatrix},  
	\begin{bmatrix} 0 & \tfrac 1\tau M+A & 0 & B^T \\
								- B & -B(\frac 1\tau M + A)^{-1} B^T & 0 & 0 \\
									0 & 0 & 0 & 0 \\
									0 & 0 & 0 & I
	\end{bmatrix}
	\bigr).
\end{equation*}
If we define $S_A := B(\frac 1\tau M^{-1} + A)^{-1} B^T$, move the second row and column to the left and bottom, respectively, and rescale certain rows and columns, then we find that the pair is equivalent to
\begin{equation*}
	\bigl(\mathcal E, \mathcal A\bigr) \sim 
	\bigl( 
	\begin{bmatrix} I & 0 & 0 & 0\\
									* & I & 0 & *_1  \\
									0 & 0 & I & *_2 \\
									0 & 0 & 0 & 0 

	\end{bmatrix},  
	\begin{bmatrix} 0 & 0 & * & *_3 \\
									0 & 0 & 0 & 0  \\
									0 & 0 & * & 0  \\
								  *_4 & 0 & 0 & -S_A 
	\end{bmatrix}
	\bigr),
\end{equation*}
where the $*$'s stand for unspecified but possibly nonzero entries. Since $S_A$ is invertible for sufficiently small $\tau$, we find that this matrix pair has the very same structure as the one of the projection scheme (see Appendix~\ref{sec_appdx_kronind_prjctn}) and, thus, is of index-1.
 \end{document}